\documentclass{amsart}

\usepackage{amsmath,rotating,graphics}
\usepackage{amsfonts}
\usepackage{amscd}
\usepackage{url}
\usepackage{amssymb,stmaryrd}
\usepackage[dvipsnames]{xcolor}
\usepackage[all]{xy}
\input xypic
\input xypic

\newcommand\restr[2]{{
  \left.\kern-\nulldelimiterspace 
  #1 
  \vphantom{\big|} 
  \right|_{#2} 
  }}
\newtheorem{defn}{Definition}[section]
\newtheorem{thm}[defn]{Theorem}

\newtheorem{prop}[defn]{Proposition}
\newtheorem{lemma}[defn]{Lemma}

\newtheorem{cor}[defn]{Corollary}

\newtheorem{rem}[defn]{Remark}
\newtheorem{schol}[defn]{Scholium}

\newtheorem{ass}[defn]{Hypothesis}
\newtheorem{ques}[defn]{Question}
\def\qed{\hfill \ensuremath{\Box}}

\newcommand{\lm}{\ensuremath{\longrightarrow}}

\newcommand{\st}{\,|\,}

\newcommand{\eps}{\varepsilon}

\DeclareMathOperator{\Hom}{\mbox{Hom}}
\DeclareMathOperator{\grhom}{\underline{\mbox{Hom}}}
\DeclareMathOperator{\grend}{\underline{\mbox{End}}}
\DeclareMathOperator{\sMCM}{\underline{MCM}}
\DeclareMathOperator{\MCM}{MCM}

\DeclareMathOperator{\Aut}{\mbox{Aut}}

\DeclareMathOperator{\Log}{\mbox{Log}}

\DeclareMathOperator{\im}{\mbox{im}\,}
\DeclareMathOperator{\spec}{\mbox{Spec}\,}

\DeclareMathOperator{\Mod}{\mbox{Mod}}
\DeclareMathOperator{\irr}{\mbox{irr}}

\DeclareMathOperator{\McK}{McK}

\DeclareMathOperator{\PGL}{PGL}
\DeclareMathOperator{\GL}{GL}
\DeclareMathOperator{\SL}{SL}

\DeclareMathOperator{\coker}{\mbox{coker}\,}

\DeclareMathOperator{\s}{\sigma}
\DeclareMathOperator{\Lam}{\Lambda}

\DeclareMathOperator{\z}{\zeta}

\DeclareMathOperator{\Z}{\mathbb{Z}}

\begin{document}

\date{\today}

\title{Low dimensional orders of finite representation type}

\author{Daniel Chan}
\address{University of New South Wales}
\email{danielc@unsw.edu.au}

\author{Colin Ingalls}
\address{
School of Mathematics and Statistics,
Carleton University, 
Ottawa, ON K1S 5B6, 
Canada}
\email{ColinIngalls@cunet.carleton.ca}

\thanks{The first author was supported by the Australian Research Council, Discovery Project Grant DP0880143. The second author was supported by an NSERC Discovery Grant.}

\begin{abstract}
In this paper, we study noncommutative surface singularities arising from orders. The singularities we study are mild in the sense that they have finite representation type or, equivalently, are log terminal in the sense of the Mori minimal model program for orders~\cite{CI}. These were classified independently by Artin (in terms of ramification data) and Reiten-Van den Bergh (in terms of their AR-quivers). The first main goal of this paper is to connect these two classifications, by going through the finite subgroups $G \subset \GL_2$, explicitly computing $H^2(G,k^*)$, and then matching these up with Artin's list of ramification data and Reiten-Van den Bergh's AR-quivers. This provides a semi-independent proof of their classifications and extends the study of canonical orders in \cite{CHI} to the case of log terminal orders. A secondary goal of this paper is to study noncommutative analogues of plane curves which arise as follows. Let $B = k_{\z}\llbracket x,y \rrbracket$ be the skew power series ring where $\zeta$ is a root of unity, or more generally a terminal order over a complete local ring. We consider rings of the form $A = B/(f)$ where $f \in Z(B)$ which we interpret to be the ring of functions on a noncommutative plane curve. We classify those noncommutative plane curves which are of finite representation type and compute their AR-quivers.
\end{abstract}

\maketitle

{\footnotesize\tableofcontents}

Throughout, we work over an algebraically closed base field $k$ of characteristic zero.

\section{Introduction}  \label{sintro}

The study of surface singularities is a classical topic related to group theory and Lie theory. It has prompted fascinating theories such as the McKay correspondence which connect these different disciplines, but it seems the correct setting for this requires one to broaden our study to include noncommutative surface singularities which arise as 2-dimensional orders. These have now already been studied by many authors including Artin, who classified maximal orders of finite representation type \cite{A} in terms of ramification data, and Reiten-Van den Bergh, who more generally classified normal orders of finite representation type in terms of their AR-quivers \cite{RV}.  These orders were also studied in \cite{DrRo94} and \cite{LVV}. In joint work with Hacking \cite[Theorem~5.4]{CHI}, we showed that these are precisely the log terminal orders on surfaces as defined in the minimal model program for orders \cite{CI}. The first goal of this paper is to reprove and connect Artin and Reiten-Van den Bergh's classifications by linking both to a classification of pairs $(G,\eta)$ where $G$ is a finite subgroup of $\GL_2$ and $\eta \in H^2(G,k^*)$. This, in  particular, extends the detailed study of canonical orders in \cite{CHI} to the case of log terminal orders. The second goal of this paper is to study noncommutative analogues of plane curve singularities $B/(f)$ where $B$ is a terminal $k\llbracket u,v \rrbracket$-order such as the skew power series ring $k_{\z}\llbracket x,y \rrbracket = k\langle \langle x,y \rangle\rangle /(yx - \z xy)$ for $\z$ some root of unity and  $0 \neq f \in Z(B)$. We classify those of finite representation type and compute their AR-quivers. 

The classification of normal orders $A$ of finite representation type is up to reflexive Morita equivalence. As observed in \cite{RV}, Section~5.3, any such $A$ is reflexive Morita equivalent to a crossed product order $k\llbracket x,y \rrbracket*_{\eta} G$ where $G$ is a finite subgroup $G$ of $\GL_2$ and $\eta \in H^2(G,k^*)$. Our first goal is to classify such pairs $(G,\eta)$. The classification of the subgroups $G$ is classical and is described in several places.  We will follow Coxeter's book \cite{Cox}, which is based on Du Val \cite{DuV}.  The nonabelian groups are classified by quadruples $(ab,a:G_1,G_2)$ where $G_1$ is a finite subgroup of $\SL_2$,  and $G_2$ is a normal subgroup with $G_1/G_2$ cyclic of order $b$ and $a$ is a positive integer (details are in Section~\ref{sCoxclass}). Roughly speaking, $G$ is the quotient of the fibre product $\mu_{ab}\times_{\mu_b} G_1$ by $\langle (-1,-1)\rangle$. We first compute $H^2(G,k^*)$ in terms of the data $a,b,G_1,G_2$. This is done by using Bergman's diamond lemma \cite{Ber} and spectral sequences to analyse group extensions. The theory for this is set up in Section~\ref{sextdiamond}. The classification does involve some case-by-case analysis, but where possible, we have tried to avoid this. 

Having classified pairs $(G,\eta)$, we can connect this classification with Reiten-Van den Bergh's by computing the AR-quivers of the skew group algebras $k\llbracket x,y \rrbracket*_{\eta} G$ which are just the McKay quivers of the pair $(G,\eta)$. Recall that the AR-quivers appearing in \cite{RV} all have the form $\Z\Delta/\langle \rho\rangle$ where $\Delta$ is an extended Dynkin quiver and $\rho$ is an automorphism of the associated translation quiver $\Z\Delta$. In Sections~\ref{sARsplit} and \ref{sARgeneral}, we compute formulas for $\Delta,\rho$ in terms of the data $a,b,G_1,G_2$. Note that \cite{RV} computes the possible relations on the AR quivers and we do not compute the relations on the resulting McKay quivers.  Fortunately there is usually no choice for the relations, so they will be uniquely determined, but this not always the case.

The next goal is to connect the classification of pairs $(G,\eta)$ to Artin's classification of ramification data. In principle, this is easy since the ramification data associated to $k\llbracket x,y \rrbracket*_{\eta} G$  is given by $\eta$ and the ramification data of the quotient morphism $\spec k\llbracket x,y \rrbracket \lm \spec k\llbracket x,y \rrbracket^G$. Of particular interest is the case where the centre $k\llbracket x,y \rrbracket^G$ is smooth. The corresponding groups $G$ are those generated by pseudo-reflections as classified by Shephard-Todd \cite{ST}. In Section~\ref{slogterm}, we give an independent classification of the possible ramification data when the centre is smooth by computing discrepancies and using the fact that normal orders of finite representation type are precisely the log terminal ones. 

We end the paper by studying the aforementioned noncommutative plane curves $B/(f)$ where $B$ is a terminal order and $f \neq 0$ is a non-unit in $Z(B)$. To determine which ones have finite representation type and compute their AR-quivers, we reformulate Kn\"{o}rrer's version of Eisenbud's theory of matrix factorisation. Here, the stable category of matrix factorisations is re-interpreted as the stable category $\sMCM \Lambda$ of maximal Cohen-Macaulay modules over the skew group algebra $\Lambda = B[\sqrt{f}] * \Z/2$ where $B[\sqrt{f}] = B[z]/(z^2-f)$ and the generator of $\Z/2$ acts by $z \mapsto -z$. The category of maximal Cohen-Macaulay $B/(f)$-modules is equivalent to $\sMCM \Lambda$ and the point is that $\Lambda$ is a 2-dimensional order so questions about finite representation type and AR-quivers can now be answered easily with our results on log terminal orders. In particular, we see this is one of the many instances where it is more natural to consider the noncommutative setting. Considering matrix factorisations really takes us out of the commutative setting into the noncommutative setting of orders.  Similar constructions are used in \cite{BFI}.

Note that specifying a noncommutative curve by giving the terminal order $B$ and defining equation $f$ is tantamount to giving the ramification data of the auxiliary order $\Lambda$. Hence relating Artin's classification to Reiten-Van den Bergh's classification is key to determining AR-quivers of noncommutative plane curves.

\noindent
\textbf{Notation:} We let $\mu_m=\mu_m(k)$ be the group of $m^{th}$ roots of unity in $k$.  We will let $\zeta_m$ denote a choice of generator of $\mu_m$. If we need also introduce $\zeta_{mn}$ we will assume that it has been compatibly chosen, namely, so that $\zeta_{mn}^n = \zeta_m$. For example, if $k = \mathbb{C}$, we may as well let $\zeta_m = e^{2 \pi i/m}$. Since in any particular question, we will only ever deal with a finite number of roots of unity, this will not be a problem.

\noindent
\textbf{Acknowledgements:} We would like to thank Anthony Henderson for his help regarding automorphisms of Dynkin quivers.   The second author would like to thank Osamu Iyama for useful conversations.  We would also like to thank an anonymous referee for many helpful suggestions.

\section{Background}  \label{sbackground}

The first half of this paper is a continuation of the study of 2-dimensional complete local orders of finite representation type studied, amongst others, by Artin \cite{A86}, Reiten-Van den Bergh \cite{RV}, \cite{LVV}, \cite{DrRo94}, and, more recently, in joint work with Hacking in \cite{CHI}. This section serves to remind the reader of the most relevant material from these papers and place the results of this paper in context. The reader can consult the first four sections of \cite{CHI} for further details. 

\subsection{McKay quivers}  \label{ssMcKay}


Recall that given any group $G$ and finite dimensional representation $V$, we can form its McKay quiver $\McK(G,V)$ whose vertices are the irreducible isomorphism classes of $G$-modules and the number of arrows from a representation $\lambda$ to $\lambda'$ is given by the multiplicity of $\lambda'$ in $\lambda \otimes_k V$. 

To cater for orders with non-trivial Brauer class, we will need a mild generalisation of this concept. We first need to recall some basic results about $k^*$-universal central extensions as introduced by Schur, details of which can be found in \cite[\S~11E]{CR}. Consider the finite abelian cohomology group $H^2(G,k^*)$ and its dual $H^2(G,k^*)^{\vee}$. Let the top row, in the following morphism of central extensions, be a $k^*$-universal central extension.
\[ \label{eSchurcover}
\begin{CD}
1 @>>> H^2(G,k^*)^{\vee} @>>> E_G @>>> G @>>> 1 \\
@. @V{\eta}VV @VVV @VV{\text{id}}V @. \\
 1 @>>> k^* @>>> \widetilde{G} @>>> G @>>> 1
\end{CD}
\]
Given any $\eta \in H^2(G,k^*)$, pushing out by the morphism $\eta: H^2(G,k^*)^{\vee} \lm k^*$ induces the central extension of $G$ by $k^*$ in the bottom row and universality means that all such central extensions arise this way.  The group $E_G$ is a Schur cover of $G$. Note that the map $\eta$ factors through its image $\mu_e < k^*$, so $\widetilde{G}$ is also induced by a central extension of $G$ by $\mu_e$. In other words, $\eta \in H^2(G,k^*)$ is induced by an element of $H^2(G, \mu_e)$. 

We look now at the McKay quiver of $(E_G,V)$. Since $H^2(G,k^*)^{\vee}$ is central in $E_G$, it acts via a character $\eta \in H^2(G,k^*)$ on any irreducible $E_G$-module $\lambda$. We define the {\em McKay quiver of $(G,\eta)$} to be the full subquiver $\McK_\eta(G,V) \subseteq \McK(\widetilde{G},V)$ consisting of these vertices. These vertices of course can be identified with the set of isomorphism classes of irreducible $\widetilde{G}$-modules where $k^*$ acts by scalar multiplication.  We will often  abbreviate $\McK_\eta(G,V)$ to $\McK_\eta(G)$. 

\begin{prop}  Let $G$ be a finite group acting linearly on a finite dimensional vector space $V$ and $E_G$ be its Schur cover as above. Then the McKay quiver of $E_G$ partitions as 
$$\McK(E_G,V) = \bigcup^{\centerdot}_{\eta \in H^2(G,k^*)} \McK_\eta(G,V).$$
Furthermore, if $V$ is a faithful representation, then the $\McK_{\eta}(G,V)$ are connected components.
  \end{prop}
\noindent\textbf{Proof.}
First, note that the $\McK_{\eta}(G,V)$ are disjoint for if $\lambda$ is a vertex in $\McK_{\eta}(G,V)$, then $H^2(G,k^*)^{\vee}$ acts on $\lambda \otimes V$ via $\eta$ too. It thus suffices to show that $\McK_{\eta}(G,V)$ is connected under the assumption that $V$ is faithful. To this end, let $\lambda, \lambda'$ be representations in $\McK_{\eta}(G)$. We need to show that there exists a non-negative integer $i$ such that 
$$\Hom_{E_G}(\lambda',\lambda \otimes V^{\otimes i}) = \Hom_{E_G}(\lambda' \otimes \lambda^* , V^{\otimes i})\neq 0.$$
Now $\lambda' \otimes \lambda^*$ is a $G$-module, so we are done by the Brauer-Burnside Theorem \cite[Theorem 4.3]{Isaacs}. 
\qed

\subsection{Crossed Product Algebras}  \label{scrossproduct}
Let $G$ be a finite group and $\eta \in H^2(G,k^*)$. Consider a two dimensional representation $\rho \colon G \lm \GL_2$. Associated to this data is the crossed product algebra $k\llbracket x,y \rrbracket*_{\eta} G$, which will play a key role in our work.  The elements of $k\llbracket x,y \rrbracket*_\eta G$ are uniquely given by elements of the form $\sum_{g \in G} f_g g$ where the $f_g\in k\llbracket x,y \rrbracket$.  We can present $\eta$ by a cocycle  $\eta : G \times G \to \mu_e$ which we also call $\eta$.  Then multiplication satisfies the relation
$$ f_ggf_hh = f_gg(f_h)\eta(g,h)gh.$$

We consider first the case where $\eta$ is trivial in which case we drop the subscript $\eta$ from the notation. Let $K$ be the field of fractions of $k\llbracket x,y \rrbracket$. Maschke's theorem tells us that $K*G$ is a separable $K^G$-algebra. Let $H = \ker \rho$ and note that $G$ acts on $K*G$ and $Z(kH)$ by conjugation.
\begin{lemma}  \label{lgeneralcross}
The primitive central idempotents of $K*G$ are precisely those of the form $\eps = \eps_1 + \ldots + \eps_r$ where $\{\eps_1, \ldots, \eps_r\} $ is a $G$-orbit of central primitive idempotents of $kH$. If $E$ denotes the set of these idempotents, then
$$k\llbracket x,y \rrbracket *G \simeq \prod_{\eps \in E} \eps k\llbracket x,y \rrbracket *G.$$
\end{lemma}
\noindent\textbf{Proof.}
Note that
$$ Z(K*G) = (Z(KH))^G = (K \otimes_k Z(kH))^G$$ 
from which the description of primitive central idempotents follows. The last assertion follows from the fact that all the idempotents lie also in $k\llbracket x,y \rrbracket*G$.
\qed

Suppose now that $G < \GL_2$ and $\eta$ is no longer trivial. As in \S~\ref{ssMcKay}, $\eta$ can be represented by some central extension $\widetilde{G}$ of $G$ by $\mu_e$ for some positive integer $e>1$. The quotient map defines a representation $\rho\colon \tilde{G} \lm \GL_2$ and we may apply the lemma to this situation. Here $\ker \rho = \mu_e$ and we may pick the primitive central idempotent $\eps \in k\mu_e$ corresponding to the standard representation of $\mu_e$ on $k$. Then one easily sees that  $k\llbracket x,y \rrbracket *_{\eta} G = \eps k\llbracket x,y \rrbracket*\widetilde{G}$. Its centre is $R= k\llbracket x,y \rrbracket^G$. 

Let $R$ be more generally, a complete local Cohen-Macaulay commutative ring and $\Lam$ be an $R$-algebra which is finitely generated as an $R$-module and has centre $R$. We say that a finitely generated $\Lam$-module $M$ is {\em maximal Cohen-Macaulay} if it is maximal Cohen-Macaulay as an $R$-module, that is $\text{depth}_R M = \dim R$. We say that $\Lam$ has {\em finite representation type} if (up to isomorphism) the number of indecomposable maximal Cohen-Macaulay modules is finite.

We will initially be interested in the case where $R$ is a 2-dimensional (commutative) complete local normal domain $R$ with field of fractions $K$. In this case, the maximal Cohen-Macaulay modules are precisely the reflexive modules. Recall that an {\em $R$-order} is an $R$-algebra $A$ such that i) $A$ is a finitely generated torsion-free $R$-module and ii) $K \otimes_R A$ is a central simple algebra. We say $A$ is {\em normal} if it is reflexive as an $R$-module and further, satisfies the following analogue of Serre's regularity in codimension one condition.
\begin{enumerate}
 \item For any height one prime $P \triangleleft R$, the localisation $A_P$ is hereditary and,
 \item the Jacobson radical $\text{rad} A_P$ is a principal left and right ideal.   
\end{enumerate}
Dropping condition~ii) above leads to the notion of tame orders as studied by Reiten-Van den Bergh \cite{RV}. 
The crossed product orders above provide the most important example of normal orders for this paper. The reason is the following result. Auslander-Reiten (AR) quivers are defined in Chapter 5 of~\cite{Yoshino} or~\cite{RV}.
\begin{thm} \cite[Theorem~5.6]{RV} 
Every tame order $A$ of finite representation type is reflexive Morita equivalent to a crossed product algebra of the form $k\llbracket x,y \rrbracket*_{\eta} G$ some pair $(G,\eta)$ as above. In particular, $A$ and $k\llbracket x,y \rrbracket*_{\eta} G$ have the same AR-quivers. 
\end{thm}
Recall that given a representation $\rho \colon G \lm \GL_2$ of a finite group, the module category of $k\llbracket x,y \rrbracket* G$ is equivalent to the category of quasi-coherent sheaves on the stack $[\spec k\llbracket x,y \rrbracket /G]$. 
Conversely, by \cite[Lemma~2.2.3]{AV}, any smooth complete local Deligne-Mumford stack of dimension two has this form. By complete local, we mean here that the coarse moduli scheme is the spectrum of a complete local ring. Now Maschke's theorem tells us that $k\llbracket x,y \rrbracket*G$ has global dimension two. Hence from Lemma~\ref{lgeneralcross}, we know that it is a product of tame orders of finite representation type. The theory of smooth Deligne-Mumford stacks is thus intimately related to the theory of tame orders of finite representation type. 

\subsection{Ramification Data}  \label{ssramdata}
The background for this material can be read in~\cite{Ch,CI,GrIn}.
The theory of ramification for orders is most naturally presented in the context of normal orders. Let $A$ be a normal $R$-order.
The {\em ramification} or {\em discriminant locus} of $A$ is the codimension one locus of the closed subset $D \subset \spec R$ where $A$ is not Azumaya.
Given any ramification curve $C$, that is, irreducible component of $D$ say corresponding to the prime $P \triangleleft R$, we can define the {\em ramification index of $A$ at $C$} to be
$$e_{C} := \dim_{R/P} Z(A_P/\text{rad}\, A_P).$$
The ramification indices of $A = k\llbracket x,y \rrbracket *_{\eta} G$ are the same as that of the morphism of commutative schemes $\spec k\llbracket x,y \rrbracket \lm \spec k\llbracket x,y \rrbracket^G$, and these in turn are the orders of the inertia groups of $G$ acting on $\spec k\llbracket x,y \rrbracket$. Indeed, from the theory of orders, we know that $e_C$ is the integer such that 
$(\text{rad}\, A_P)^{e_C} = PA_P$ where $P \triangleleft k\llbracket x,y \rrbracket^G$ is the prime corresponding to $C$. Let $Q\triangleleft k\llbracket x,y \rrbracket_P$ be semiprime such that $Q^e = Pk\llbracket x,y \rrbracket_P$. Now $k\llbracket x,y \rrbracket_P/Q *_{\eta} G$ is a separable algebra so $\text{rad}\, A_P = QA_P$ and we conclude $e=e_C$ as desired.

Ramification data give a natural way of associating a log surface to a normal order $A$. We define the boundary divisor 
$$\Delta_A = \sum_C \left(1 - \frac{1}{e_C}\right)C$$
and define $\Log(A):=(\spec Z(A), \Delta_A)$ to be the {\em log surface associated to $A$}. Note that $\Log(A)$ only depends on the reflexive Morita equivalence class of $A$. This is because ramification data is determined by codimension one data, and reflexive Morita equivalence corresponds to a Morita equivalence on the punctured spectrum of $R$. 

In \cite{CI}, we introduced a version of Mori's minimal model program  for orders on surfaces. We will not reproduce the definition of discrepancy for orders (see \cite[Definition~3.13]{CI}) from which follow the definitions of terminal, canonical and log terminal orders. It suffices for now to note the following result.
\begin{thm}\cite[Proposition~3.15]{CI} \cite[Theorem~5.4]{CHI}
The following are equivalent conditions on a normal order $A$.
\begin{enumerate}
 \item $A$ has finite representation type.
 \item $A$ is log terminal.
 \item $(\spec Z(A), \Delta_A)$ is a log terminal log surface. 
\end{enumerate}
\end{thm}
In light of the above result, it is natural to ask the following question.
  \begin{ques}
    Let $R$ be a complete local $k$-algebra of dimension two and let $\Delta$
    be a divisor on $\spec R$ with standard coefficients, i.e. of the form $(1-1/e)$ with $e \in \Z, e \geq 1$.  If we suppose that $(\spec R,\Delta)$ is log terminal
    does this imply that there exists a finite subgroup $G\subseteq \GL_2$ such that $R \simeq k[[x,y]]^G$
    and $\Delta$ is the ramification divisor of $k[[x,y]] \to R$?
  \end{ques}
  
    \section{The classification of subgroups of $\GL_2$} \label{sCoxclass}

In this section, we recall Coxeter's description of finite subgroups of $\GL_2$ as well as Shephard-Todd's descriptions of finite groups generated by pseudo-reflections. We start with the former and consider $G < \GL_2$. Scalar multiplication of matrices induces a double cover of groups $\pi: k^* \times \SL_2 \lm \GL_2$  so it suffices to classify the double covers $\widetilde{G}:=\pi^{-1}(G)$. Now the finite subgroups of $k^*$
are all of the form $\mu_m$ for some $m$ whilst the finite subgroups of $\SL_2$ are 
classified by the Dynkin diagrams $A_n,D_n,E_{6,7,8}$.  Let $p,q$ be the projections $\SL_2 \times k^* \to \SL_2$ and $\SL_2 \times k^* \to k^*$.  Let $G_1 =p(\widetilde{G}),$ $G_2 = \ker(\restr{q}{\tilde{G}})$, and $\mu_{m} = q(\widetilde{G})$.  Then $\widetilde{G}$ is a subgroup of $G_1 \times \mu_{m}$ and so by Goursat's Lemma \cite{G} we 
 are thus reduced to classifying diagrams of the form 
$$ \begin{CD}
  1   & @>>> & G_2      & @=   & G_2 \\
 @VVV &      & @VVV     &      & @VVV \\
\mu_a & @>>> &  \widetilde{G}       & @>>> & G_1 \\
 @|   &      & @VVV     &      & @VVV \\
\mu_a & @>>> & \mu_{ab} & @>>> & \mu_b 
\end{CD} $$
where $m=ab$ and the 3 term sequences are short exact sequences and 
the lower right square is a fibre product. In particular, 
$\widetilde{G} \simeq \mu_{ab} \times_{\mu_b} G_1$ and $G \simeq \widetilde{G}/\langle (-1,-1) \rangle$.

The possibilities for $G_1, G_2$ above are fairly limited since 
$G_1$ is a finite subgroup of $\SL_2$ and $G_2$ is 
a normal subgroup with a cyclic quotient so $[G_1,G_1] \leq
G_2 \leq G_1$. The table below lists all such 
possibilities and is taken from \cite[Section~10.1]{Cox}.
$$\begin{array}{l|l|l|l|l}
  type        &  b   & a    &  G_1     &  G_2 \\
  \hline
A_{m,n}^c      &   b   &  a        & A_{db-1} & A_{b-1}     \\
B_n^m       & 2    & 2m   & D_{n+2}  & A_{2n-1} \\
D_n^m       & 1    & 2m   & D_n      & D_n     \\
CD_n^m      & 2    & 2m   & D_{2n}   & D_{n+1} \\
BT_n^{2m+1} & 4    & 2m+1 & D_{2n+3} & A_{2n}  \\
E_6^m       & 1    & 2m   & E_6      & E_6     \\
E_7^m       & 1    & 2m   & E_7      & E_7     \\
E_8^m       & 1    & 2m   & E_8      & E_8     \\
F_{41}^m    & 2    & 2m   & E_7      & E_6     \\
G_{21}^m    & 3    & 2m   & E_6      & D_4     
\end{array}$$
In the non-abelian cases, that is all cases except type $A$, we have $b \leq 4$ so the homomorphism $G_1 \lm \mu_b$ is defined up to inverse and we may denote the corresponding subgroup $G < \GL_2$ by $G = (Z_1,Z_2;G_1,G_2)$ where $Z_1 = ab$ and $Z_2 = a$. We note the following fact which can be checked case-by-case. 
\begin{rem}  \label{rsplitmodminus}
We can find $g_1 \in G_1$ generating $G_1/G_2$ such that $g_1^b = \pm I$ where $I$ is the identity matrix. 
\end{rem}

In the abelian case, there are a lot of choices for $G_1 \lm \mu_b$ so we classify these directly instead as follows.  First note that we may diagonalise $G$ and so assume 
$G \subset k^* \times k^*$. Taking quotients we obtain an exact
sequence
$$ 1 \rightarrow G \rightarrow k^* \times k^* \stackrel{M}{\rightarrow} k^* \times k^*\rightarrow 1 $$
where $M$ is an integer matrix with non-zero determinant.  We may change $M$ via
left multiplication by $\GL_2(\Z)$ to get it in the form 
$$M =\begin{pmatrix} m & 0 \\
                     -c & n  
\end{pmatrix},$$
where $m,n \geq 0$ and $0 \leq c < m$.
In this case $G$ is generated by the matrices
$$\begin{pmatrix} \z_m & 0 \\
                     0 & 1  
\end{pmatrix}\quad \quad \begin{pmatrix} \z_{mn}^c & 0 \\
                     0 & \z_{n}  
\end{pmatrix}$$
We will denote this group by $A_{m,n}^c$. Its dual group $G^{\vee} \simeq \coker(\Z^2 \xrightarrow{M^T} \Z^2)$ which has order $\det(M) = mn$.  We write $p \vee q$ for the positive greatest common divisor of $p$ and $q$. Since $G$ is non-canonically isomorphic to $G^{\vee}$, the theory of invariant factors shows that 
$$A_{m,n}^c \simeq \Z/(mn/(m \vee n \vee c)) \oplus \Z/(m \vee n \vee c).$$ 
Its subgroup generated by pseudo-reflections is isomorphic to $\Z/m \,\, \oplus \,\, \Z/(n \vee c)$. Also, $\det(G) = \langle \zeta_{mn}^n, \zeta_{mn}^{m+c} \rangle = \langle \zeta_{mn}^{n\vee(m+c)}\rangle$ so the order of the cyclic group $SG := G \cap \SL_2$ is $n \vee (m+c)$. We have not computed the relationship between the parameters $m,n,c$ with the parameters $a,b,d$ since we will have no need for it. 

There are many other notations for various finite subgroups of $\GL_2$ and we recall the ones we will use here. Given positive integers $p \geq q \geq r$ we follow Coxeter \cite{Cox} in defining 
$$ \langle p,q,r \rangle := \langle \alpha, \beta, \gamma, \zeta \st 
\alpha^p = \beta^q = \gamma^r = \alpha \beta \gamma = \zeta, \zeta^2 = 1 \rangle.$$
In particular, $\langle 3, 3, 2\rangle$, $\langle 4, 3, 2 \rangle$ $\langle 5, 3, 2 \rangle$ are the symmetry groups of the Platonic solids. 

Shephard-Todd classify the finite unitary groups generated by pseudo-reflections in \cite{ST}. We will also need both their classification and notation. Any such group $G$ is either primitive or imprimitive. The primitive ones form 19 exceptional groups numbered 4-22 in \cite[Section~4]{ST}, and which we will denote by ST4-ST22. The imprimitive ones are the $G(np,p,2), n,p \in \mathbb{Z}_{>0}$ defined as follows. Let $H$ be the subgroup of $\mu_{np} \times \mu_{np}$ consisting of $(\zeta, \xi)$ such that $\zeta\xi \in \mu_n$. Then $G(np,p,2)$ is the group generated by the diagonal image of $H$ in $\GL_2$ and the involution
$$
\begin{pmatrix}
0 & 1 \\ 1 & 0
\end{pmatrix}.
$$

\section{Group extensions via the diamond lemma}  \label{sextdiamond} 

Given a finite subgroup $G$ of $\GL_2$, we would like to give a nice description of the cohomology group $H^2(G,k^*)$ and the corresponding central extensions. Our criterion for ``nice'', is one which removes as much case-by-case analysis as possible, and facilitates the computation of the McKay quiver (see Section~\ref{sARsplit}). One natural approach is to use the fact that $G$ fits in an exact sequence of the form
\begin{equation}
1 \lm H \lm G \lm \overline{G} \lm 1
\label{eG}
\end{equation}
where $\overline{G}$ is a cyclic group of order say $n$ and $H$ has simpler group cohomology than $G$. Then one can appeal to the Hochschild-Serre spectral sequence
$$ E^{pq}_2 = H^p(\overline{G},H^q(H,k^*)) \Longrightarrow H^{p+q}(G,k^*).$$ 
To make this work, we have to properly understand this sequence and especially, the differentials $d_2,d_3$. We devote this section towards this goal. Note first that the spectral sequence gives an exact sequence 
\begin{equation}
 0 \lm  H^1(\overline{G},H^{\vee})^{d_2=0} \xrightarrow{\iota} H^2(G,k^*) \xrightarrow{\rho}  H^0(\overline{G},H^2(H,k^*))^{d_2=0} \xrightarrow{d_3} H^3(\overline{G},k^*) 
\label{ess}
\end{equation}
where $ H^1(\overline{G},H^{\vee})^{d_2=0}$ denotes 
$\ker (d_2: H^1(\overline{G},H^{\vee}) \lm H^3(\overline{G},H^0(H,k^*))$
and similarly for the second last term. From Section~\ref{sCoxclass}, we may choose $H$ to be a subgroup of $\SL_2$ so $H^2(H,k^*) = 0$ and a natural question is how to write the central extension of $G$ given a 1-cocycle representing an element of  $ H^1(\overline{G},H^{\vee})^{d_2=0}$.

Now any (not necessarily central) extension of $G$ is certainly also an extension of $\overline{G}$ which motivates the need for 
\begin{lemma}  \label{ldiamond}  
Let $\widetilde{H}$ be a finite group and $\overline{G}$ a cyclic group of order $n$. Let $\phi \in \Aut \widetilde{H}$ and $h_0 \in \widetilde{H}$ be such that 
\begin{enumerate}
\item $\phi(h_0) = h_0$ and
\item  $\phi^n$ is conjugation by $h_0$.
\end{enumerate} 
We define the group 
$$ \widetilde{G} := \langle \widetilde{H},\widetilde{g} \,\st\, \widetilde{g}^n = h_0, \widetilde{g} h= \phi(h) \widetilde{g}\ \text{for}\ h \in \widetilde{H}\rangle.$$
Then there is an exact sequence
$$1 \lm \widetilde{H} \lm \widetilde{G} \lm \overline{G} \lm 1$$
where $\widetilde{H}$ embeds naturally in $\widetilde{G}$ and $\widetilde{g}$ maps to a generator of $\overline{G}$. 

Conversely, any group extension of $\overline{G}$ by $\widetilde{H}$ as above comes from such a construction.
\end{lemma}
\noindent\textbf{Proof.} The converse direction is easy to see for given $\widetilde{G}$ we may arbitrarily lift a generator of $\overline{G}$ to $\widetilde{g} \in \widetilde{G}$ and let $h_0 = \widetilde{g}^n$ and $\phi$ be conjugation by $\widetilde{g}$. 

We hence prove the forward direction by showing every element of $\widetilde{G}$ can be written uniquely in the form $\widetilde{g}^jh$ for some $h \in \widetilde{H}, j = 0, \ldots, n-1$. This follows from Bergman's diamond lemma on checking overlaps as we shall now verify. Given $h,h' \in \widetilde{H}$, the overlap check for $\widetilde{g}^nh$ is
$$ \phi^n(h)h_0 = \phi^n(h)\widetilde{g}^n = \widetilde{g}^n h = h_0 h $$
so is verified precisely when condition ii) above holds. Similarly the overlap check for $\widetilde{g}^{n+1}$ corresponds to condition i) and the overlap check for $\widetilde{g}hh'$ corresponds to the fact that $\phi$ is a group homomorphism. $\qed$

We begin by fixing an exact sequence of finite groups of the form (\ref{eG}) where $\overline{G}$ is cyclic of order $n$ but $H$ is arbitrary. We fix a generator $\overline{g} \in \overline{G}$ as well as a lift to $g \in G$. Hence $g^n =: h_0 \in H$. As is commonly done, given a $\overline{G}$-module $M$, we will compute the cohomology groups $H^i(\overline{G}, M)$ by the cohomology of the complex
$$ 0 \lm M \xrightarrow{D} M \xrightarrow{N} M \xrightarrow{D} M \xrightarrow{N} \ldots $$
where $D = 1 - \overline{g}$ and $ N = 1 + \overline{g} + \ldots + \overline{g}^{n-1}$. In particular, when we refer to cocycles, we mean with respect to this complex. 

We now define our candidates for the differentials in the Hochschild-Serre spectral sequence. Let $\chi \in H^{\vee}$ be a 1-cocycle representing an element $[\chi]$ of $H^1(\overline{G}, H^{\vee})$. Let 
\begin{equation}
d_2: H^1(\overline{G}, H^{\vee}) \lm H^3(\overline{G}, k^*) = \mu_n: [\chi] \mapsto \chi(h_0) 
\label{ed2}  
\end{equation}
which is easily checked to be well-defined. Consider now a central extension
$$\eta: 1 \lm k^* \lm \widetilde{H} \xrightarrow{\pi} H \lm 1$$
representing an element of $H^0(\overline{G},H^2(H,k^*))$. Since $\eta$ is $\overline{G}$-invariant, conjugation by $g$ on $H$ lifts to an automorphism $\phi \in \Aut \widetilde{H}$ which fixes $k^*$. Note that any character $\chi \in H^{\vee}$ induces an automorphism $\widetilde{h} \mapsto \chi(\pi(\widetilde{h}))\widetilde{h}$ of $\widetilde{H}$ fixing $k^*$, and that $\phi$ is determined only up to an element of $H^{\vee}$. We lift $h_0$ to an element $\widetilde{h_0} \in \widetilde{H}$. Now conjugation by $\widetilde{h_0}$ induces an automorphism $\phi_0 \in \Aut \widetilde{H}$ which is independent of the lift and differs from $\phi^n$ by a character. We may hence define
$$d_2:H^0(\overline{G},H^2(H,k^*)) \lm H^2(\overline{G},H^{\vee}):\eta \mapsto [\phi^n \circ \phi^{-1}_0] .$$
Now $g$ commutes with $h_0$ so $\phi(\widetilde{h_0})$ and $\widetilde{h_0}$ differ by a scalar in $k^*$ which is furthermore in $\mu_n=H^3(\overline{G},k^*)$ if $d_2(\eta) = 0$. We may thus define
$$ d_3:H^0(\overline{G},H^2(H,k^*))^{d_2=0} \lm H^3(\overline{G},k^*): \eta \mapsto \phi(\widetilde{h_0})\widetilde{ h_0}^{-1}.$$

Rather than verifying that the above are the differentials in the Hochschild-Serre spectral sequence, we will verify exactness of the sequence (\ref{ess}) which is what we will actually use.

\begin{thm}  \label{tHSseq}  
Consider an exact sequence $1 \lm H \lm G \lm \overline{G} \lm 1$ where $\overline{G}$ is cyclic of order $n$ and let $d_2,d_3$ be the maps defined above. Then there is an exact sequence
$$  0 \lm  H^1(\overline{G},H^{\vee})^{d_2=0} \xrightarrow{\iota} H^2(G,k^*) \xrightarrow{\rho}  H^0(\overline{G},H^2(H,k^*))^{d_2=0} \xrightarrow{d_3} H^3(\overline{G},k^*) $$
\end{thm}
\noindent\textbf{Proof.} We will only sketch the proof here as it follows by the usual arguments from Lemma~\ref{ldiamond}. We continue the notation above. Let $\widetilde{G}$ be a central extension of $G$ by $k^*$ which represents an element of 
$\ker(\rho: H^2(G,k^*) \lm H^2(H,k^*))$. Then $\widetilde{G}$ splits when restricted to $H$ so there is an exact sequence 
$$ 1 \lm k^* \times H \lm \widetilde{G} \lm \overline{G} \lm 1 .$$
We apply Lemma~\ref{ldiamond} with $\widetilde{H} = k^* \times H$ to show $\ker \rho =  H^1(\overline{G},H^{\vee})^{d_2=0}$ as follows. We may lift $g$ to $\widetilde{g} \in \widetilde{G}$ so that $\widetilde{g}^n = (1,g^n)=:h_0$.  Conjugation by $\widetilde{g}$ induces an automorphism $\phi \in \Aut k^* \times H$ which must have the form $\phi(\alpha,h) = (\chi(h)\alpha,ghg^{-1})$ for some character $\chi \in H^{\vee}$. Condition~ii) of Lemma~\ref{ldiamond}, namely, that $\phi^n$ is conjugation by $h_0$ corresponds to the fact that $\chi$ is a 1-cocycle representing an element of $H^1(\overline{G},H^{\vee})$. Condition~i) of Lemma~\ref{ldiamond} corresponds to the fact that $d_2([\chi]) = 0$. It follows fairly readily now that  $\ker \rho =  H^1(\overline{G},H^{\vee})^{d_2=0}$.

We now compute $\im \rho$ using Lemma~\ref{ldiamond}. We start with a central extension $\widetilde{H}$ of $H$ by $k^*$ corresponding to $\eta \in H^0(\overline{G}, H^2(H,k^*))$ and consider when it lifts to an extension of $G$ and so in particular gives an exact sequence 
$$ 1 \lm \widetilde{H} \lm \widetilde{G} \lm \overline{G} \lm 1 .$$
We first lift $g^n \in H$ to an element $\widetilde{h_0} \in \widetilde{H}$. Since $\widetilde{G}$ is also a central extension of $G$, Lemma~\ref{ldiamond} tells us that our lift of $\widetilde{H}$ to some $\widetilde{G}$ corresponds to some automorphism $\phi \in \Aut \widetilde{H}$ such that i) $\phi^n$ is conjugation by $\widetilde{h_0}$ and ii) $\phi(\widetilde{h_0}) = \widetilde{h_0}$. Now condition~i) corresponds precisely to the condition $d_2(\eta) =0$ and condition~ii) corresponds to the condition $d_3(\eta) = 0$. $\qed$

For future reference, we record the next result which follows from the proof of the theorem and Lemma~\ref{ldiamond}.

\begin{schol}  \label{sgroupext}  
Consider an element of $H^1(\overline{G},H^{\vee})^{d_2=0}$ represented by the 1-cocycle $\chi \in H^{\vee}$. The group extension corresponding to $\iota([\chi]) \in H^2(G,k^*)$ is 
$$\widetilde{G} := \langle k^* \times H, \widetilde{g} \st \widetilde{g}^n = (1,g^n), \widetilde{g} (\alpha,h) = (\chi(h)\alpha, ghg^{-1}) \widetilde{g},\ \text{for}\ \alpha \in k^*, h \in H \rangle.$$
\end{schol}

\section{McKay quivers \`a la Auslander-Reiten}  \label{sARsplit}  

McKay quivers of small subgroups of $\GL_2$ were computed in \cite[Section~2]{AR}. In this section, we reformulate their setup and pave the way for computing McKay quivers of general pairs $(G,\eta)$ in the next section. The relevance of McKay quivers for us is the following result, which is well known when $\eta$ is trivial.

\begin{prop}\cite[Proposition~4.5]{CHI}
Given a finite subgroup $G < \GL_2$ and $\eta \in H^2(G,k^*)$, the AR-quiver of $k\llbracket x,y \rrbracket*_{\eta} G$ is the McKay quiver of the pair $(G,\eta)$. If $W$ is an irreducible representation of $\widetilde{G}$ corresponding to a vertex of $\McK_\eta(G)$, then the corresponding vertex of the AR-quiver is the $k\llbracket u,v \rrbracket*_{\eta} G$-module $k\llbracket u,v \rrbracket \otimes_k W$. 
\end{prop}

To avoid confusion when writing elements of the skew group ring, we will let $\hat{\zeta}$ be the scalar matrix with entries $\zeta$. 

Given quivers $Q_1,Q_2$, we consider the product quiver $Q_1 \times Q_2$ whose vertices are $(v_1,v_2)$ where $v_1,v_2$ are vertices of $Q_1,Q_2$ and the arrows are $(v_1,v_2) \xrightarrow{(e_1,e_2)} (w_1,w_2)$ for arrows $v_1 \xrightarrow{e_1} w_1$ in $Q_1$ and $v_2 \xrightarrow{e_2}w_2$ in $Q_2$. We also let $\Z/a$ denote the quiver with vertices the elements of $\Z/a$ and arrows $j \rightarrow j+1$ for $j \in \Z/a\Z$. We can think of $\Z$ as the McKay quiver of $k^* = \GL_1$. Similarly, for $a>0$, $\Z/a$ is the McKay quiver of $\mu_a < \GL_1$.

We fix $\widehat{-1} \in H < \SL_2$ with McKay quiver $\Delta_H$, which is the double of an extended Dynkin quiver $\Delta$. Note that $\widehat{-1}$ lies in all finite subgroups of $\SL_2$ except $A_{2n}$.  To pick an orientation of arrows on $\Delta$, we call the {\em even} vertices of $\Delta_H$, those where $\widehat{-1}$ acts trivially and the other ones the {\em odd} vertices. Then $\Delta_H$ is bipartite with every arrow going between an even and an odd vertex. We may hence pick the orientation on $\Delta$ so that all arrows go from even to odd vertices. We similarly need to consider the {\em even} part of $\McK(k^* \times H) = \Z \times \Delta_H$ which is defined to be the full subquiver of $\McK(k^* \times H)$ on those vertices where $(-1,\widehat{-1}) \in k^* \times H$ acts trivially.
We write $k^*H$ for the subgroup of $\GL_2$ generated by $k^*$ and $H$. Hence as in \cite[Lemma~6]{AR}, we have
$$ \McK(k^* H) = (\Z \times \Delta_H)^{ev}  = \Z\Delta$$
where the superscript $ev$ denotes the even part, and $\Z\Delta$ is defined as in \cite[\S VII.4,p.~250]{ARS}. Analogously, we note that for $m \in \Z_+$ we have $\mu_{2m}^{\vee} \stackrel{\text{can}}{\simeq} \Z/2m$ so with corresponding notation we have $\McK(\mu_{2m} H) = (\Z/2m  \times \Delta_H)^{ev}$.

It was shown in \cite{RV} that AR-quivers of log terminal orders can be described as quotients of $\Z\Delta$ by an automorphism. It is best to describe these automorphisms in terms of $\Z \times \Delta_H$ rather than $\Z\Delta$. They can be gotten by the product of an automorphism of $\Z$ with an automorphism of $\Delta_H$. The only automorphisms of the quiver $\Z$ come from addition by $n \in \Z$ which we will denote by $[+n]$. 
There are two natural ways to obtain automorphisms of $\Delta_H$. Given a one dimensional representation $\chi$ of $H$, tensoring by $\chi$ permutes the irreducible representations and so induces an automorphism of $\Delta_H$. We will call such an automorphism of $\Delta_H$ the {\em character automorphism} associated to $\chi$, or sometimes somewhat ambiguously, associated to $\ker \chi$.  Another automorphism comes from considering a finite subgroup of $G<\GL_2$ which contains $H$ as a normal subgroup. For any $g \in G$, the outer automorphism of $H$ defined by conjugation by $g$ induces an automorphism of $\Delta_H$. If $g$ generates the group $G/H$ then we will say this automorphism is an {\em outer automorphism associated to $G$}. When $G/H$ is order two, then this is uniquely determined so we may speak of the outer automorphism associated to $G$. It is a good exercise in group theory to show that these character and outer automorphisms generate all automorphisms of the quiver $\Delta_H$. 

We now consider the problem of relating McKay quivers of a pair with nice a ``subpair''. Consider firstly, a finite subgroup $G<\GL_2$ and a cyclic quotient $G/\widetilde{H} =: \overline{G} \simeq \Z/n$. 
Let $\Lambda:= k\widetilde{H}$ and consider a  crossed product algebra  $\widetilde{\Lambda} := k *_{\eta}G$. Suppose that we can express $\widetilde{\Lambda} = \Lambda *_{\gamma} \overline{G}$ for some cocycle $\gamma$ with values in $Z(\Lambda)^*$. This occurs if $\gamma$ is trivial on $\widetilde{H}$ and we can find $g \in G$ with $g^n \in Z(\widetilde{H})$ as occurs for example in Remark~\ref{rsplitmodminus}. 

Let $\irr \widetilde{H}$ and $\irr (G,\eta)$ denote the set of isomorphism classes of irreducible $\Lambda$ and $\widetilde{\Lambda}$-modules. We first note that $\overline{G}$ acts on $\irr \widetilde{H}$ as follows. Given any automorphism $\phi$ of $\Lambda$ and irreducible $\Lambda$-module $M$, we obtain a new irreducible $\Lambda$-module $M_{\phi}$ which has the same underlying addition as $M$, but $\lambda \in \Lambda$ acts by $\phi(\lambda)$. Inner automorphisms act trivially on $\irr \widetilde{H}$ so $\overline{G}$ acts via outer automorphisms induced by conjugation. Similarly, $\overline{G}^{\vee}$ acts on $\widetilde{\Lambda}$ and hence $\irr (G,\eta)$. Indeed, $\widetilde{\Lambda}$ is naturally $\overline{G}$-graded so given $\overline{g} \in \overline{G}$ the character $\chi \in \overline{G}^{\vee}$ acts on the degree $\overline{g}$-graded component by $\chi(\overline{g})$. 

The key to computing McKay quivers is the following standard lemma. 
\begin{lemma}  \label{lindres}
With the above notation,
  \begin{enumerate}
\item if $\overline{G}$ acts freely on $\irr \widetilde{H}$ then $\McK_\eta(G) = \McK(\widetilde{H})/\overline{G}$.
\item if $\overline{G}^{\vee}$ acts freely on $\irr (G,\eta)$ then $\McK(\widetilde{H}) = \McK_\eta(G)/\overline{G}^{\vee}$.
\end{enumerate}
\end{lemma}
\noindent\textbf{Proof.} This exercise in induction and restriction follows as in the proof of Proposition~1.8 of \cite{RR}. $\qed$

The following is a uniform statement of \cite[Proposition~7]{AR}, giving the McKay quiver of a group. We include their proof as it will be referred to later. Below as usual, we use Coxeter's notation for subgroups of $\GL_2$ as recalled in Section~\ref{sCoxclass}. 

\begin{prop} \label{pARgroup} 
Let $G = (ab,a;G_1,G_2)$ be a non-abelian subgroup of $\GL_2$ and $\rho$ be the character automorphism of $\Delta_{G_1}$ associated to the character $G_1 \lm \mu_b$. Then $\McK(G) = (\Z\times \Delta_{G_1})^{ev}/\langle [+a] \times\rho \rangle$. \end{prop}
\noindent\textbf{Proof.} When $b=1$, we have that $\rho$ is trivial and $G=\mu_{a}G_1$. Its McKay quiver, as already noted, is $(\Z/a \times \Delta_{G_1})^{ev}$ so the proposition holds in this case. The general case follows from this $b=1$ case and Lemma~\ref{lindres} applied to $\Lambda= kG \hookrightarrow \widetilde{\Lambda} = k(\mu_{ab}G_1)$ noting that $\mu_b^{\vee}$ does indeed act freely on $\irr \mu_{ab}G_1$. $\qed$

\section{McKay quivers in the general case}  \label{sARgeneral}  

In this section, we classify pairs $(G,\eta)$ and compute the corresponding McKay quivers. This gives
the AR-quivers of all log terminal orders and, in particular, gives a group theoretic proof of the fact that they have the form $\Z\Delta/\langle \rho \rangle$ as proved by Reiten-Van den Bergh \cite{RV}. 

We assume to begin with that $G$ is a non-abelian finite subgroup of $\GL_2$ and that $\eta \in H^2(G,k^*)$. As noted in Section~\ref{ssMcKay}, we know this comes from an extension $\widetilde{G}$ of $G$ by $\mu_e$  for sufficiently large and divisible $e$ . Recall also that $k*_{\eta} G$ is the direct factor of $k\widetilde{G}$ whose  module category is the  subcategory where $\mu_e<\widetilde{G}$ acts as scalar multiplication. Below we let $\zeta_m$ denote an appropriate primitive $m$-th root of unity as explained in the Notation at the end of the Introduction. 

\begin{prop}  \label{pAReasybrauer}   
Let $G = (2m,2m;G_1,G_1)$. Then $H^2(G,k^*)\simeq \Hom(\mu_m,(G_1/\langle \widehat{-1} \rangle)^{\vee})$. Let $\eta \in  \Hom(\mu_m,(G_1/\langle \widehat{-1} \rangle)^{\vee})$ and $\chi = \eta(\zeta_m)$ which we consider as a character of $G_1$ of order $e$ say. Let $H = \ker \chi$ and $g_1 \in G_1$ be an element with $\chi(g_1) = \zeta_e$ and $g_1^e = \widehat{-1}$. Then $\McK_\eta(G) = (\Z\times\Delta_H)^{ev}/\langle [+2m/e] \times \rho\rangle$ where $\rho$ is the outer automorphism of $\Delta_H$ induced by conjugation by $g_1$.
\end{prop}
\noindent\textbf{Proof.} We have an exact sequence of the form
$$ 1 \lm G_1 \lm G \lm \mu_m \lm 1 $$
so we may use Theorem~\ref{tHSseq} to identify $H^2(G,k^*) = H^1(\mu_m, G_1^{\vee})^{d_2 = 0}$. Now $\hat{\zeta}_{2m} \in G$ is a lift of the generator $\zeta_m \in \mu_m$ so $\mu_m$ acts trivially on $G_1^{\vee}$. Moreover, given the description of $d_2$ in Section~\ref{sextdiamond}, the elements of $ H^1(\mu_m, G_1^{\vee})^{d_2 = 0}$ are precisely the characters $\chi \in G_1^{\vee}$ of order $m$ such that $\chi(\hat{\zeta}_{2m}^m = \widehat{-1}) = 1$. This shows that $H^2(G,k^*)\simeq \Hom(\mu_m,(G_1/\langle \widehat{-1} \rangle)^{\vee})$.

The group extension corresponding to $\eta$ is by Scholium~\ref{sgroupext},
$$ \widetilde{G} = \langle \mu_e \times G_1 , \widetilde{\zeta}\st  \widetilde{\zeta}^m = (1,\widehat{-1}), \widetilde{\zeta}(\alpha,g) = (\alpha \chi(g),g)\widetilde{\zeta} \ \text{for}\ \alpha \in \mu_e, g \in G_1 \rangle .$$
Note that this contains the subgroup $\mu_{2m}H$ on identifying $H$ with $1 \times H$ and $\hat{\zeta}_{2m}$ with $\widetilde{\zeta}$. We apply Lemma~\ref{lindres} to $\Lambda = k(\mu_{2m}H), \widetilde{\Lambda} = k*_{\eta} G$. We need to determine how conjugation by $g_1\in \widetilde{\Lambda}$ operates on $\Lambda$ and hence $\irr \mu_{2m}H$. Now in $\widetilde{\Lambda}$ we have 
$$ g_1\widetilde{\zeta} g_1^{-1} = \chi(g_1) \widetilde{\zeta} = \zeta_{2m}^{2m/e} \widetilde{\zeta} $$
so $g_1$ acts freely and the McKay quiver is as stated in the proposition. $\qed$

We give the McKay quivers for the remaining cases where $G$ is non-abelian in the following 

\begin{prop}  \label{pARhardbrauer}  
Let $G = (ab,a;G_1,G_2)$ be a non-abelian group with $G_2 \neq G_1$. Then one of the following occurs.
\begin{enumerate}
\item $(G_1,G_2) = (D_{2n+3},A_{2n}), (E_7,E_6)$ or $(E_6,D_4)$. In this case $H^2(G,k^*) = 0$.
\item $G = (4m,2m;D_{n+2},A_{2n-1})$. Then $H^2(G,k^*) \simeq \Z/2$ if $m$ is odd and is zero otherwise. If $\eta$ represents the non-zero element we have $\McK_\eta(G) = (\Z \times \Delta_{A_{2n-1}})^{ev}/\langle [+m] \times \rho_1\rho_2 \rangle$ where $\rho_1$ is a character automorphism associated to $1$ and $\rho_2$ is the outer automorphism associated to $D_{n+2}$.
\item $G = (4m,2m;D_{2n},D_{n+1})$ and $m+n$ is even. Then $H^2(G,k^*) = 0$.
\item $G = (4m,2m;D_{2n},D_{n+1})$, $n$ is even and $m$ is odd. Then $H^2(G,k^*) \simeq \Z/2$ and if $\eta$ represents the non-zero element we have $\McK_\eta(G) = (\Z \times \Delta_{D_{n+1}})^{ev}/\langle [+m] \times \rho_1 \rho_2 \rangle$ where $\rho_1$ is a character automorphism associated to $A_{n-2}$ and $\rho_2$ is the outer automorphism associated to $D_{2n}$.  
\item $G = (4m,2m;D_{2n},D_{n+1})$, $n$ is odd and $m$ is even.  Then $H^2(G,k^*) \simeq \Z/2$ and if $\eta$ represents the non-zero element we have $\McK_\eta(G) = (\Z \times \Delta_{D_{n+1}})^{ev}/\langle [+m] \times \rho_1 \rho_2 \rangle$ where $\rho_1$ is the character automorphism associated to $D_{(n+3)/2}$ and $\rho_2$ is the outer automorphism associated to $D_{2n}$. 
\end{enumerate}
\end{prop}
\noindent\textbf{Proof.} The possibilities for $G$ were recalled in Section~\ref{sCoxclass} while the computation of $H^2(G,k^*)$ follows directly from Theorem~\ref{tHSseq} and is left to the reader.  

We prove cases iv) and v) only. Case ii) is similar and easier. Let $\s,\tau$ be generators for $D_{2n}$ with relations $\s^{2n-2} = \tau^2 = -1, \tau \s = \s^{-1}\tau$ so that $D_{n+1} = \langle \s^2,\tau \rangle$. Suppose first that $n$ is even and $m$ is odd. We consider the exact sequence
\begin{equation}  
 1 \lm \widetilde{H}:= \mu_{2m}D_{n+1} \lm G \lm \overline{G} \lm 1
\label{eDTseq} 
\end{equation}
where $\overline{G}\simeq \Z/2\Z$. We may lift the generator of $\overline{G}$ to $g = \hat{\z}_{4m}\s\tau \in G$. Incidentally, one can prove in this case that $H^2(G,k^*) \simeq \Z/2\Z$ by showing $\ker (d_3:H^0(\overline{G},H^2(H,k^*)) \lm H^3(\overline{G},k^*)) = 0$. We seek a non-trivial element in $H^1(\overline{G}, \widetilde{H}^{\vee})^{d_2 = 0}$ which must then give $\eta$. Note that $\widetilde{H}^{\vee}$ is an index two subgroup of $\mu_{2m}^{\vee} \times D_{n+1}^{\vee}$. In fact, if $\chi_1$ denotes the canonical generator of $\mu_{2m}^{\vee}$ and $\chi_2:D_{n+1} \lm k^*: \tau \mapsto \z_4, \s^2 \mapsto -1$, then $\widetilde{H}^{\vee} = \langle \chi_1\chi_2, \chi_1^2,\chi_2^2 \rangle$. Now conjugation by $g$ maps $\chi_1 \mapsto \chi_1, \chi_2 \mapsto \chi_2^{-1}$ so one computes that 
$$ H^1(\overline{G}, \widetilde{H}^{\vee}) = 
\frac{\langle \chi_1^m\chi_2,\chi_2^2\rangle}{\langle \chi_2^2\rangle}.$$
Thus $\chi := \chi_1^m\chi_2$ gives a non-trivial cohomology class and furthermore, $\chi(g^2) = \chi(-\hat{\z}_{2m}) = 1$ so $d_2([\chi]) = 0$. The group extension corresponding to $\eta$ is thus
$$\widetilde{G} = \langle k^* \times \widetilde{H}, \widetilde{g} \st
\widetilde{g}^2 = (1,-\hat{\z}_{2m}), \widetilde{g}(\xi,h) = (\xi\chi(h),ghg^{-1})\widetilde{g} \rangle.$$
We apply Lemma~\ref{lindres}i) to $k\widetilde{H} \lm k*_{\eta} G$ to obtain iv). Indeed, $\chi_1^m$ corresponds to the automorphism $[+m]$, $\chi_2$ corresponds to the character automorphism $\rho_1$ and conjugation by $g$ corresponds to the outer automorphism $\rho_2$. 

We now prove v) and assume that $n$ is odd and $m$ is even. In this case, we do not use the exact sequence~(\ref{eDTseq}) since the $E^{0,2}_{\infty}$ term does not vanish. We consider instead the exact sequence
$$1 \lm \widetilde{H} \lm G \lm \overline{G} \lm 1$$
where $\overline{G} \simeq \Z/2$ and 
$$\widetilde{H} := \langle \hat{\z}_{4m}\s\tau,\s^2 \rangle \simeq (4m,2m; D_{n+1}, A_{2n-3}).$$
Note that $\tau \in G$ lifts the generator of $\overline{G}$. We seek a non-trivial element of $H^1(\overline{G}, \widetilde{H}^{\vee})^{d_2 = 0}$. Now $\widetilde{H}^{\vee} = \langle \chi_1,\chi_2 \rangle$ where
$$ \chi_1(\hat{\z}_{4m} \s\tau) = \z_{2m}, \ \chi_1(\s^2) = 1, \ \chi_2(\hat{\z}_{4m} \s\tau) = 1, \ \chi_2(\s^2) = -1 .$$
As in case iv), one finds that $\chi = \chi_1^{m/2}\chi_2$ gives a non-trivial cohomology class corresponding to $\eta$. The corresponding group extension has the form 
$$\widetilde{G} := \langle k^* \times \widetilde{H}, \widetilde{\tau} \st
\widetilde{\tau}^2 = (1,\widehat{-1}), \widetilde{\tau} (\xi,h) = (\xi \chi(h), \tau h \tau^{-1}) \widetilde{\tau} \rangle.$$ 
Thus conjugation by $\widetilde{\tau}$ on $k\widetilde{H}$ maps $\hat{\z}_{4m}\s\tau \mapsto \z_4 \hat{\z}_{4m}\tau(\s\tau)\tau^{-1}, \s^2 \mapsto - \tau \s^2 \tau^{-1}$.
We apply Lemma~\ref{lindres}i) to $k\widetilde{H} \lm k*_{\eta} G$ to write $\McK_\eta(G)$ as a quotient of $\McK(\widetilde{H})$ which in turn is a quotient of $\McK(\mu_{4m}D_{n+1})  = (\Z/4m\Z \times \Delta_{D_{n+1}})^{ev}$. To obtain the description of the McKay quiver in v), we lift conjugation by $\widetilde{\tau}$ on $k\widetilde{H}$ to the (order 4) automorphism of $k \mu_{4m}D_{n+1}$ defined by
$$ \hat{\z}_{4m} \mapsto \z_4 \hat{\z}_{4m}, \ \s \tau \mapsto \tau (\s\tau)\tau^{-1}, \ 
\s^2 \mapsto - \tau\s^2 \tau^{-1} .$$
This induces the automorphism $[+m] \times \rho_1\rho_2$ of $\McK(\mu_{4m}D_{n+1})$ given in v), which we note also induces a free action on $\McK(k\widetilde{H})$. $\qed$

Finally, we consider now the case where $G$ is abelian. So up to conjugacy, $G$ is a diagonal subgroup of $\GL_2$ and hence a subgroup of $k^* \times k^*$. In fact, it must be of the form $\mu_{ab}\times_{\mu_a} \mu_{ac}$ for some positive integers $a,b,c$ and choice of surjections $\mu_{ab} \lm \mu_a$ and $\mu_{ac} \lm \mu_a$. 
Let $\Z \oplus \Z$ denote the McKay quiver of $k^* \times k^*$ which has vertices the elements of $\Z \oplus \Z$ and arrows those of the form $(n,n') \rightarrow (n+1,n')$ or $(n,n') \rightarrow (n,n'+1)$. Similarly, we denote the McKay quiver of $\mu_b \times \mu_c$ by $\Z/b \oplus \Z/c$. To describe the various McKay quivers that occur, we recall that any subgroup $L$ of $\Z/b \oplus \Z/c$ acts freely on the quiver $\Z/b \oplus \Z/c$ so we may form the quotient quiver $(\Z/b \oplus \Z/c) /L$. 

\begin{prop} \label{pARabelian}  
Let $G = \mu_{ab} \times_{\mu_a} \mu_{ac}$. Then $H^2(G,k^*) \simeq \mu_{d_0}^{\vee}$ where $d_0 = b\vee c$. Furthermore, given any $\eta \in H^2(G,k^*)$, the McKay quiver of $(G,\eta)$ has the form $(\Z \oplus \Z) / L$ for some rank two subgroup $L < \Z \oplus \Z$. 
\end{prop}
\noindent\textbf{Proof.} We consider first the case when $\eta = 0$. If $a=1$ then $\McK(G) = \mu_{b}^{\vee} \oplus \mu_c^{\vee}$. For general $a$, we note that there exists an exact sequence
$$ 1 \lm G \lm \mu_{ab} \times \mu_{ac} \xrightarrow{\phi} \mu_a \lm 1 $$
where $\phi$ is concocted from the surjections $\mu_{ab} \lm \mu_a, \mu_{ac} \lm \mu_a$. Dualising this sequence, we see that $\McK(G) = (\mu_{ab}^{\vee} \oplus \mu_{ac}^{\vee})/\mu_a^{\vee}$.

We now compute $H^2(G,k^*)$ using Theorem~\ref{tHSseq} applied to the exact sequence
$$ 1 \lm \mu_b \lm G \lm \mu_{ac} \lm 1 .$$
Pick $g = (\zeta_{ab},\zeta_0) \in G$ and note $\zeta_0$ generates $\mu_{ac}$. Then 
$$H^2(G,k^*) = H^1(\mu_{ac},\mu_b^{\vee})^{d_2=0} = \Hom(\mu_{ac},\mu_b^{\vee})^{d_2=0} .$$
Elements of the latter correspond to characters $\chi:\mu_b \lm k^*$ such that $\chi^{ac} = 1$ and $\chi(\zeta_b^c) = 1$. Now $\zeta_b^c$ generates the subgroup $\mu_{b/d_0}$ so the two conditions on $\chi$ just amount to the fact that $\chi$ factors through $\mu_{d_0}$. Hence $H^2(G,k^*) \simeq \mu_{d_0}^{\vee}$. Given such a $\chi$ representing $\eta \in H^2(G,k^*)$, the corresponding group extension is given by 
$$ \widetilde{G} = \langle k^* \times \mu_b, \widetilde{g}\st \widetilde{g}^{ac} = (1,\zeta_b^c), \widetilde{g}(\xi,\zeta) = (\xi\chi(\zeta),\zeta)\widetilde{g}, \text{for}\ \xi \in k^*, \zeta \in \mu_b \rangle .$$
Suppose that $\chi$ has order $d$ so that the subgroup $\widetilde{H}:= \langle 1 \times \mu_{b/d}, \widetilde{g} \rangle \simeq \mu_{ab} \times_{\mu_{ad}} \mu_{ac}$. We may apply Lemma~\ref{lindres}i) to $k\widetilde{H} \hookrightarrow k*_{\eta} G$. Indeed, from the relations for $\widetilde{G}$, we see that conjugation by $(1,\zeta_b)$ on $k\widetilde{H}$ is an automorphism of the form $\widetilde{h} \mapsto \lambda(\widetilde{h}) \widetilde{h}$ for some character $\lambda \in \widetilde{H}^{\vee}$. Then noting that we have already proved $\McK(\widetilde{H}) = \widetilde{H}^{\vee}$ we see that $\McK_\eta(G) = \widetilde{H}^{\vee}/\langle \lambda \rangle \simeq \Z \oplus \Z /L$ for some appropriate $L$.  $\qed$

\textbf{Remark:} 
\begin{enumerate}
\item The proof of Proposition~\ref{pARabelian} shows how to compute the McKay quiver in this case. We have not given a description of it since there is no succinct way to write it down uniformly. 
\item The computations for $H^2(G,k^*)$ in Propositions~\ref{pAReasybrauer}, \ref{pARhardbrauer} and \ref{pARabelian} are tabulated in the appendix.
\item Sections~\ref{sARsplit} and \ref{sARgeneral} complete the computation of all AR-quivers of log terminal orders in terms of pairs $(G,\eta)$. Reiten-Van den Bergh \cite{RV} classify all these orders in terms of the AR-quiver and mesh-type relations on the path algebra of this quiver. When $\Delta$ is not type $A$, there is, up to isomorphism, only one choice for these relations. In type $A$, the relation is either uniquely determined, or depends on a parameter $\z$ which is a root of unity in $k$. It would be an interesting exercise to determine $\z$ as a function of the $(G,\eta)$ in this case. 
\end{enumerate}

\section{Log terminal $k\llbracket u,v \rrbracket$-orders} \label{slogterm} 

In this section, we classify log terminal orders with smooth centre via ramification data. We do this directly by computing discrepancies as opposed to appealing to Shephard-Todd's classification of finite groups generated by pseudo-reflections. We then elucidate the correspondence between finite subgroups of $\GL_2$ and the possible ramification data of log terminal orders. 

We work more generally with a log surface $(\spec k\llbracket u,v \rrbracket,\Delta)$ where $\Delta$ has {\em standard coefficients} in the sense that $\Delta = \sum_i b_i C_i$ where each coefficient $b_i$ has the form $1-1/e_i$ for some positive integer $e_i$ which we call the {\em ramification index}. Note that the log surface associated to a normal $R$-order always has standard coefficients. We will usually write the ramification indices $e_1,e_2, \ldots,e_m$ with multiplicity equal to the multiplicity of the corresponding curves $C_i$ in $\Delta$. 

Before giving the following necessary condition for an order to be log terminal, we need to recall some definitions. Given an effective divisor $\Delta = \sum_i b_i C_i$ with $C_i$ irreducible, the {\em support of } $\Delta$ is $\text{Supp}\ \Delta = \cup_i C_i$. Of particular interest, are {\em simple singularities,} originally classified in \cite{Arnold}, or for an exposition see for example, \cite[Definition~8.1]{Yoshino}. 
Following \cite[Chapter~II, \S~8]{BHPV}, these can be characterised as curves $C \subset \spec k\llbracket u,v \rrbracket$ which have multiplicity $\leq 3$ and such that for any iterated blowup at closed points $\pi: Z' \lm \spec k\llbracket u,v \rrbracket$, the {\em reduced total transform} $\text{Supp}\ \pi^*C$ consists only of points of multiplicity $\leq 3$.  We say that positive integers $a,b,c$ form a {\em Platonic triple} if $1/a + 1/b + 1/c > 1$ so this holds in particular, in the degenerate case where one of the integers is 1.

Below we use the basic theory and defintions of the minimal model program as found in ~\cite{KM}.  In particular, {\em klt} stands for Kawamata log terminal as defined in \cite[Defn.~2.34]{KM}.
\begin{lemma}  \label{lmult}  
Let $(\spec k\llbracket u,v \rrbracket,\Delta)$ be a log surface with standard coefficients and $e_1, \ldots,e_m$ be the ramification indices written with multiplicity. Let $\pi:Z_1 \lm \spec k\llbracket u,v \rrbracket$ be the blowup at the closed point and $E$ the exceptional curve.
\begin{enumerate}
\item The discrepancy of $(\spec k\llbracket u,v \rrbracket,\Delta)$ along $E$ is greater than -1 if and only if $m \leq 3$ and when $m=3$, we have that $e_1,e_2,e_3$ form a Platonic triple. 
\item If $(\spec k\llbracket u,v \rrbracket,\Delta)$ is klt then $\text{Supp}\ \Delta$ is a simple singularity.
\end{enumerate}
In particular, the ramification locus of a log terminal $k\llbracket u,v \rrbracket$-order is a simple singularity.
\end{lemma}
\noindent\textbf{Proof.} The discrepancy of $\Log(\Lam)$ along $E$ is 
$$a_1 = 1- \sum_{i=1}^m \left(1-\frac{1}{e_i}\right).$$
Hence $a_1 > -1$ if and only if 
\begin{equation}
\sum \frac{1}{e_i} > m-2.
\label{emult} 
\end{equation}
Now the left hand sum is at most $m/2$ so $m\leq 3$. When $m=3$ we also see that the inequality in (1) reduces to $\sum 1/e_i > 1$ so the ramification indices form a Platonic triple. Part~i) follows. 

For part ii), note first that points of multiplicity $\leq 2$ are simple so we may assume that $\text{Supp}\ \Delta$ has a triple point. If $\pi^{-1}_*$ denotes strict transform, then \cite[Lemma~2.30]{KM} tells us that $(Z_1,\pi^*\Delta)=(Z_1, \pi^{-1}_*\Delta + (2-\sum 1/e_i)E)$ 
is klt so the same is true of $(Z_1, \pi^{-1}_*\Delta + (1- 1/e_i)E)$ for any $i$ by \cite[Corollary~2.35(1)]{KM}. Part~ii) now follows from part~i) by induction. $\qed$

We now classify log terminal orders via ramification data. We continue with the above notation. Also, our notation for simple curve singularities will follow that of~\cite[Proposition 8.5]{Yoshino}.

\begin{thm} \label{tlogterm}  
Let $\Lambda$ be a normal $k\llbracket u,v \rrbracket$-order and $C$ its ramification locus. Let $e_1,\ldots,e_m$ be the ramification indices written with multiplicity.  Then $\Lambda$ is log terminal if and only if its ramification is one of the following.
\begin{enumerate}
\item $C$ has multiplicity $\leq 1$.
\item $C$ is type $A_{2k-1}$ and $\{e_1,e_2,k\}$ form a Platonic triple.
\item $C$ is type $A_{2k}$ with ramification index $e_1$ and both  $\{e_1,e_1,k\}$ and $\{2,e_1,2k+1\}$ form Platonic triples.
\item $C$ is type $D_{2k+2}$ and if $e_3$ is the ramification index along the component which is not tangential to the other two, then $\{e_1,e_2,ke_3\}$ form a Platonic triple.
\item $C$ is type $D_{2k+3}$
  and the ramification along the singular component is $2$.
\item $C$ is type $E_6$ 
  and the ramification index is 2.
\item $C$ is type $E_7$ 
  and all ramification indices are 2. 
\item $C$ is a type $E_8$ 
  and the ramification index is 2. 
\end{enumerate}
\end{thm}
\noindent\textbf{Proof.} Case~i) gives examples of terminal orders which are log terminal by definition (see \cite{CI}). Cases~ii),iii) can be obtained from cases iv),v) as the degenerate situation where $e_3 = 1$ so we will omit their proof. 

We now prove cases iv) and v). Now $C$ is the union of a double point $C'$ with a smooth curve $C''$. In case iv), $C'$ is the union of two possibly tangential curves, whilst in in case v), it is a cusp. 
We construct a sequence of blowups 
$$ X_k \lm X_{k-1} \lm \ldots \ldots X_1 \lm \spec k\llbracket u,v \rrbracket $$
by repeatedly blowing up the double point until it becomes a smooth curve.  Let $E_j \subset X_k$ be the exceptional curve corresponding to the $j$-th blowup.

We prove case~iv) first.  The discrepancy of $\Log(\Lam)$ along  $E_j$ is 
$$a_j = j - j(1-1/e_1) - j(1-1/e_2) - (1-1/e_3) .$$
The log terminal condition is equivalent to 
$$-1 < a_j \Longleftrightarrow 1/e_1 + 1/e_2 + 1/(je_3) > 1$$
for $j = 1, \ldots, k$. This is satisfied precisely when $\{e_1,e_2,ke_3\}$ form a Platonic triple. 

To prove case~v), we will need to further blowup $X_k$ twice to obtain a log resolution. As before, let $a_j$ denote that the discrepancy of $\Log(\Lam)$ along the exceptional curve created on the $j$-th blowup. Then 
\begin{align*}
a_j  & = j - 2j(1-1/e_1) - (1-1/e_3) , \text{for}\ j \leq k\\
a_{k+1} & = (k+1) - (2k+1)(1-1/e_1) - (1-1/e_3)  \\
a_{k+2} & = (2k+2) - (4k+2)(1-1/e_1) - 2 (1-1/e_3)
\end{align*}
The smallest discrepancy is either $a_k$ or $a_{k+2}$ and $a_k> -1$ if and only if $\{e_1,e_1,ke_3\}$ form a Platonic triple while $a_{k+2} > -1$ if and only if $\{2,e_1,(2k+1)e_3\}$ form a Platonic triple.

Next, we prove case~vi). We blowup the triple point 3 times to obtain a log resolution from which we can compute the discrepancy. Let $a_1,a_2,a_3$ be the discrepancies along the 3 exceptional curves. From lemma~\ref{lmult}, we know that $a_1 > -1$ if and only if $\{e_1,e_1,e_3\}$ form a Platonic triple. Now 
$$a_2 = 2 - 3(1-1/e_1) - 2(1-1/e_3) \ \text{and} \ 
a_3 = 4-6(1-1/e_1) - 3(1-1/e_3) .$$
Hence $\Lam$ is log terminal if and only if $e_1 = e_3=2$. 

For cases~vi) and viii), we know the ramification index must be 2 by lemma~\ref{lmult}i). A calculation like those above shows that this does give a log terminal order. 

Lemma~\ref{lmult} and the classification of simple singularities now shows that there are no other possibilities for log terminal $k\llbracket u,v \rrbracket$-orders. $\qed$

Artin's classification of log terminal orders in \cite{A} proceeds via reduction to the case where the centre is smooth. Theorem~\ref{tlogterm} essentially recovers this part of his result by quite different means. 

We wish of course to match up Artin's classification in terms of ramification data to the classification of pairs $(G,\eta)$ and hence to Reiten-Van den Bergh's classification. To do this, recall from Section~\ref{ssramdata} that the ramification of $k\llbracket x,y \rrbracket *_{\eta} G$ is given by the ramification of $\spec k\llbracket x,y \rrbracket \lm \spec k\llbracket x,y \rrbracket^G \simeq k\llbracket u,v\rrbracket$. The groups giving smooth centre are those generated by pseudo-reflections classified by Shephard-Todd. Below, we first go through Shepherd-Todd's primitive groups \cite[Section~4]{ST} and record their ramification curves. These were also computed in \cite[Section~2]{B}.

$$ \begin{array}{l|l|l}
\text{ST} & \mbox{type}    & \text{ramification} \\
\hline
4  & G_{21}^1 &       u^3-v^2,3   \\
5  & E_6^3    &       u^4-v^2,3,3 \\ 
6  & G_{21}^2 &       u^6-v^2,3,2 \\
7  & E^6_6    &       v(u+v)(u-v),2,3,3 \\
8  & F_{41}^2 &       u^3-v^2,4 \\
9  & E_7^4    &       u^6-v^2,2,4 \\
10 & F_{41}^6 &       u^4-v^2,4,3 \\
11 & E_7^{12} &       v(u+v)(u-v),4,3,2 \\
12 & F_{41}^1 &       u^4-v^3,2 \\
13 & E_7^2    &       v(u^3-v^2),2,2 \\
14 & F_{41}^3 &       u^8-v^2,2,3 \\
15 & E_7^6    &       u(u^2-v)(u^2+v),2,2,3 \\
16 & E_8^5    &       u^3-v^2,5 \\
17 & E_8^{10} &       u^6-v^2,2,5\\
18 & E_8^{15} &       u^4-v^2,5,3 \\
19 & E_8^{30} &       v(u-v)(u+v),5,3,2\\
20 & E_8^3    &       u^5-v^2,3\\
21 & E_8^6    &       u^{10}-v^2,3,2\\
22 & E_8^2    &       u^5-v^3,2
\end{array}
$$

For the imprimitive groups we have the following. 
\begin{thm}
The ramification of $G(mn,n,2)$ is given in the table below. 

$$\begin{array}{l|l}
  \hline
m,n           & \emph{ramification}   \\
m=1, n \geq 1 & u^n-v^2,2 \\
m \geq 2, n\geq 1    & u(u^n-v^2),m,2 
\end{array}$$
\end{thm}

\section{Ramification of log terminal orders}

In this section, we determine the ramification of the log terminal order $k\llbracket x,y \rrbracket *_{\eta} G$ and hence match up our classification of pairs $(G,\eta)$ with Artin's classification of ramification data. 
We begin by determining an efficient way of computing the subgroup of $G = (ab,a;G_1,G_2)$ generated by pseudo-reflections. Let $\phi: G_1 \lm \mu_b$ be the homomorphism with kernel $G_2$ used to define $G$. We say $g_1 \in G_1$ {\em lifts to a pseudo-reflection} in $G$ if there is a pseudo-reflection $g \in G$ such that $g$ maps to $g_1$ under the natural projection map $G \lm G_1$. In this case, we say that $g$ is a {\em pseudo-reflection lift} of $g_1$. 

The existence of pseudo-reflection lifts is easy to analyse. Note first that we can change basis so that $g_1 = \left( \begin{smallmatrix} \xi & 0 \\ 0 & \xi^{-1}  \end{smallmatrix} \right)$. The only possible pseudo-reflection lifts are $\xi g_1$ and $\xi^{-1} g_1$. These lie in $G$ if and only if $\xi^{\pm a} = \phi(g_1)$. In particular, we have the following.

\begin{prop}  \label{preflectionlift}  
Suppose that $G$ is non-abelian so that $b \leq 4$. Then $g_1 \in G_1$ lifts to a pseudo-reflection in $G$ if and only if $\zeta^a = \phi(g_1)$ for some primitive $d$-th root of unity $\zeta\in \mu_{ab}$ where $d$ is the order of $g_1$. 
\end{prop}

As an example, we now compute the pseudo-reflection subgroup $RG$ of $G = B^m_n = (4m,2m;D_{n+2},A_{2n-1})$, as well as the corresponding quotient $G/RG$ as a subgroup of $\GL_2$. Note $|G| = 4mn$. We fix our copy of $D_{n+2}$ to be the group generated by 
$$ \sigma = 
\begin{pmatrix}
 \zeta_{2n} & 0 \\ 0 & \zeta_{2n}^{-1}
\end{pmatrix}, 
\tau = 
\begin{pmatrix}
 0 & \zeta_4 \\ \zeta_4 & 0 
\end{pmatrix}.$$
Proposition~\ref{preflectionlift} shows that the anti-diagonal elements $\sigma^j\tau$ lift to pseudo-reflections in $G$ if and only if $m$ is odd. Let us first assume that $m$ is odd so that $RG$ contains the group generated by these pseudo-reflection lifts, that is, the dihedral group generated by $\sigma$ and $\rho = \left(\begin{smallmatrix} 0 & 1 \\ 1 & 0 \end{smallmatrix}\right)$. The $\sigma^j$ which lift to pseudo-reflections are those for which $\zeta_{2n}^j \in \mu_{2m}$. Now $\mu_{2n} \cap \mu_{2m} = \mu_{2(m\vee n)}$ so the diagonal subgroup of $RG$ is precisely $\mu_{m\vee n} \times \mu_{m\vee n}$ where $m\vee n$ is the greatest common divisor of $m$ and $n$. We see immediately that $RG = G(2n, \frac{2n}{m\vee n},2)$ in the Shephard-Todd notation recalled in Section~\ref{sCoxclass}. Note that $|RG| = 4n(m \vee n)$. It will be convenient to let $\overline{m} = m/m \vee n, \overline{n} = n/m \vee n$. 

To compute $G/RG$ we first note that $k[x,y]^{RG} = k[u = x^{2n} + y^{2n}, v = (xy)^{m \vee n}]$. Also, we must have $|G/RG| = \overline{m}$. We claim that $\zeta_{2m} \sigma RG$ generates $G/RG$. Its action on $k[u,v]$ is given by $u \mapsto \zeta_{2m}^{2n} u, \ v \mapsto \zeta_{2m}^{2(m \vee n)}$. Hence $\zeta_{2m} \sigma RG$ generates the diagonal subgroup $\frac{1}{m}(n, m \vee n) = \frac{1}{\overline{m}}(\overline{n}, 1)$ which has order $\overline{m}$. This proves the claim.

We now look at the case where $m$ is even where the pseudo-reflection subgroup is the diagonal subgroup $RG = \mu_{m\vee n} \times \mu_{m\vee n}$. In this case $k[x,y]^{RG} = k[u = x^{m \vee n}, v = y^{m \vee n}]$. 
In this case, $G/RG$ is not cyclic, so we examine the action of generators for $G$ on $ku \oplus kv$.
$$ 
\zeta_{4m} \tau : \quad 
\zeta_{4m}^{m \vee n} \zeta_4^{m \vee n}\begin{pmatrix}
 0 & 1 \\ 1 & 0 
\end{pmatrix}, \quad
\sigma: \quad 
\begin{pmatrix}
 \zeta_{2n}^{m \vee n} & 0 \\ 0 & \zeta_{2n}^{-m \vee n}
\end{pmatrix}.
$$
We wish to write out this small subgroup in Coxeter notation $G/RG = (Z_1,Z_2;D_{\overline{n} + 2}, G_2)$. Considering the action of $(\zeta_{4m}\tau)^2$ shows that $Z_2 = 2\overline{m}$. 
They key question is thus, which roots of unity $\zeta$ are such that $\zeta \left(\begin{smallmatrix}
0 & \zeta_4 \\ \zeta_4 & 0 \end{smallmatrix}\right)$ lie in $G/RG$? This depends on the parity of $\overline{m}$ so we first assume that $\overline{m}$ is odd. Now $m$ is even, so $\zeta_4^{m \vee n} = \pm 1$. We thus see that $\zeta_{4m}^{m \vee n} \zeta_4^{m \vee n}\zeta_{4}^{-1}$ is trivial in $\mu_{4\overline{m}}/\mu_{2\overline{m}}$ so $\zeta$ can be any root in $\mu_{2\overline{m}}$. It follows that when $\overline{m}$ is odd, we have $G/RG = (2\overline{m},2\overline{m};D_{\overline{n}+2}, D_{\overline{n}+2})$. 
Suppose now that $\overline{m}$ is even. Now $\zeta_4 \in \mu_{2\overline{m}}$ so  $\zeta_{4m}^{m \vee n} \zeta_4^{m \vee n}\zeta_{4}^{-1}$ is non-trivial in $\mu_{4\overline{m}}/\mu_{2\overline{m}}$. It follows that now we have $G/RG = (4\overline{m},2\overline{m};D_{\overline{n}+2}, A_{2\overline{n}-1})$.

Going through the list in Section~\ref{sCoxclass} of non-abelian subgroups of $\GL_2$ and performing a similar computation gives
\begin{thm}  \label{tramofklt}
Let $G = (ab,a;G_1,G_2)$ be a non-abelian subgroup of $\GL_2$. Then the reflection subgroup $RG$ of $G$ and corresponding quotient $G/RG$ are as listed in the Appendix. 
\end{thm}
The reflection subgroups of abelian subgroups were already noted in Section~\ref{sCoxclass}.

\textbf{Remark:} In Artin's classification of ramification data of maximal orders of finite representation type, he lists the possible centres in the form $k\llbracket s,t \rrbracket^{G_0}$ by listing the groups $G_0$. Also, rather than describing the ramification loci and indices on $k\llbracket s,t \rrbracket^{G_0}$, he describes their pullback to $k\llbracket s,t \rrbracket$. Our results allow us also to describe the ramification data of $k\llbracket x,y \rrbracket *_{\eta} G$ in this format. Indeed, the group $G_0$ is just the quotient group $G/RG$ whilst the pullback of the ramification data to $k\llbracket s,t \rrbracket^{G_0}$ is the ramification associated to $RG$ listed in Section~\ref{slogterm}. Note that since Artin assumes his orders are maximal, there are fewer possibilities than what we find. For example, his group $G_0$ is always cyclic, whereas we have seen that this may not be the case when $G = B^m_n$.

\section{Kn\"orrer's method}  \label{sKnorrer}  

In this section, we review Kn\"orrer's reformulation \cite{Kn} of Eisenbud's theory of matrix factorisations in our noncommutative setting. It allows us to relate the question of finite representation type of noncommutative plane curves to a similar question about 2-dimensional orders. 

Let $R$ be a complete local commutative Cohen-Macaulay ring and $\Lam$ be an $R$-algebra which is finitely generated as an $R$-module and has centre $R$. We let $\MCM \Lam$  denote the category of maximal Cohen-Macaulay $\Lam$-modules. We will say that $\Lam$ is a {\em regular $R$-algebra} if furthermore, $R$ is regular local, $\Lam$ is free as an $R$-module and $\text{gl.dim} \Lam = \dim R$. We have the following version of the Auslander-Buchsbaum formula from \cite[Lemma 2.16]{IW}.

\begin{prop}  \label{pAB}  
Let $\Lam$ be a regular $R$-algebra. Then for any finitely generated $\Lam$-module $M$ we have 
$$ \text{depth}_R M  + \text{pd}_{\Lam} M = \dim R.$$
\end{prop}

Assume from now on that $R$ is a complete regular local ring with residue field $k$ and dimension $\geq 2$ so is isomorphic to a power series ring and that $B$ is a regular $R$-algebra. Let $0 \neq f \in R$ and consider the double covers $R[\sqrt{f}] = R[z]/(z^2 - f)$ and $B[\sqrt{f}] = B[z]/(z^2-f) = R[\sqrt{f}] \otimes_R B$.  Note that the cyclic group $G_f = \langle \s \rangle$ of order 2 acts on $R[\sqrt{f}]$ and $B[\sqrt{f}]$ by $\s: z \mapsto -z$ and fixing elements of $R$ and $B$. We may thus form the skew group ring $\Lam := B[\sqrt{f}] * G_f$. Let $\eps = \frac{1}{2}(1 - \s)\in R[\sqrt{f}]*G_f \subset \Lam$ which we note is an idempotent. Now $B[\sqrt{f}]\simeq \Lam(1-\eps)$ is a Cohen-Macaulay $\Lam$-module so we may also consider the quotient category $\sMCM \Lam = (\MCM \Lam) / \langle B[\sqrt{f}] \rangle$. We will denote homomorphism and endomorphism groups in this category with $\grhom, \grend$ respectively. The theory of matrix factorisations \`a la Eisenbud \cite[Chapter~6]{Eis} and Kn\"orrer's method \cite[Proposition~2.1]{Kn} carry over to our setting painlessly. We will use the following re-interpretation, which is also used in \cite{BFI}.

\begin{prop}  \label{pmfact}  
Let $\Lam = B[\sqrt{f}] * G_f$ be as above and $A = B/(f)$. Then there is an equivalence of categories 
$$ \grhom_{\Lam}(\Lam\eps,-): \sMCM \Lam \xrightarrow{\sim} \MCM A.$$
In particular, $A$ has finite representation type if and only if $\Lam$ does. 
\end{prop}
\noindent\textbf{Proof.} We merely sketch the proof so that the reader can easily check that the now standard arguments in \cite{Eis} and \cite{Kn} do extend to this setting.

Consider the components of the Peirce decomposition of $\Lam$ given by $\eps$,  
$$ \eps \Lam \eps = B, \eps \Lam (1-\eps) = zB, (1-\eps) \Lam \eps = zB, (1-\eps) \Lam (1-\eps) = B .$$
Note that for a $\Lam$-module $M$ we have 
$$ \grhom_{\Lam}(\Lam\eps,M) = \frac{\eps M}{\eps \Lam (1-\eps) M}  = \frac{\eps M}{z (1-\eps) M}  = \text{coker} \bigl( (1-\eps) M \xrightarrow{z} \eps M \bigr) $$
so in particular $\grend_{\Lam} \Lam \eps = \frac{\eps\Lam \eps}{z(1-\eps)\Lam\eps} = B/(f) = A$ and $\grhom_{\Lam}(\Lam \eps, M)$ is an $A$-module. Furthermore, we see now that for $M$ maximal Cohen-Macaulay, $\grhom_{\Lam}(\Lam \eps, M)$ is also maximal Cohen-Macaulay since multiplication by $z$ must be injective on $M$. 

To construct the inverse functor consider a maximal Cohen-Macaulay $A$-module $N$. By the version of the Auslander-Buchsbaum formula in Proposition~\ref{pAB}, we see that $N$ has a $B$-projective resolution of the form
$$ 0 \lm M_1 \xrightarrow{\phi} M_0 \lm N \lm 0 .$$
As in \cite[Proposition~2.1]{Kn}, the $B$-module $M_1 \oplus M_0$ can be made into a unique $\Lam$-module $M$ with $\eps M = M_0, (1-\eps)M = M_1$ and multiplication by $z$ restricts to $\phi: (1-\eps)M \lm \eps M$. Standard arguments may now be applied to show that this is indeed a well-defined functor which is inverse to $\grhom_{\Lam}(\Lam\eps,-)$.  $\qed$

Below we use the notion of normal orders recalled in Section~\ref{sbackground}. 
\begin{lemma} \label{lreduced}  
Let $B$ be a normal, regular $R$-order. If $A:=B/(f)$ has finite representation type, then $A$ has no non-zero nilpotent ideals. In particular, the subscheme $C_f:f=0\subset \spec R$ defines a reduced subscheme of $\spec R$ which does not contain any components of the ramification locus of $B$. 
\end{lemma}
\noindent\textbf{Proof.} Suppose $A$ has finite representation type. Then Auslander's Theorem \cite[Section~1, Theorem, p.~9]{Aus} ensures that $A$ is an isolated singularity, that is, for each height one prime $P \triangleleft R/(f)$, the localisation $A_P$ is regular and thus semisimple. If $A$ had a non-zero nilpotent ideal, then the same would be true of some $A_P$, a contradiction. We may thus assume that $f=0$ defines reduced subscheme of $\spec R$. Let $Q \triangleleft R$ be a prime over which $B$ ramifies say with ramification index $e$. From ramification theory, we know that $(\text{rad}\, B_Q)^e = QB_Q$ so the ideal $I = B \cap \text{rad}\, B_Q$ is nilpotent modulo $Q$. It follows that $Q$ does not divide $f$.$\qed$

In the light of the previous lemma, it makes sense to consider
\begin{ass} \label{aram}  
Let $B$ be a normal, regular $k\llbracket u,v \rrbracket$-order. Suppose further $C_f:f=0 \subset \spec k\llbracket u,v \rrbracket$ defines a reduced curve and that it does not contain any of the ramification curves of the $R$-order $B$. This ensures in particular that $\Lam:= B[\sqrt{f}]*G_f \simeq (R[\sqrt{f}]*G_f) \otimes_R B$ is a normal $k\llbracket u,v \rrbracket$-order with the combined ramification of that of $R[\sqrt{f}]/R$ and $B/R$.
\end{ass}

The regular orders $B$ of interest for us in this paper, are the terminal orders introduced in \cite{CI}. As mentioned in Section~\ref{sbackground}, these are defined via a version of discrepancy for orders. In the complete local case, these terminal orders have been classified as follows. Let $\zeta$ be a primitive $e$-th root of unity and $k_{\zeta}\llbracket x,y \rrbracket$ be the skew power series ring $k\langle \langle x,y\rangle\rangle/(yx - \zeta xy)$. It has centre $k\llbracket u,v \rrbracket$ where $u=x^e, v=y^e$. For any positive integer $n$, we define the $k\llbracket u,v \rrbracket$-order
$$
B_{e,ne} = 
\begin{pmatrix}
 k_{\zeta}\llbracket x,y \rrbracket & \ldots & \ldots & k_{\zeta}\llbracket x,y \rrbracket \\
 y k_{\zeta}\llbracket x,y \rrbracket & \ddots & \ddots & \vdots \\
\vdots & \ddots & \ddots & \vdots \\
 yk_{\zeta}\llbracket x,y \rrbracket & \ldots & y k_{\zeta}\llbracket x,y \rrbracket & k_{\zeta}\llbracket x,y \rrbracket
\end{pmatrix}
\subseteq M_n(k_{\zeta}\llbracket x,y \rrbracket)
$$
We now summarise the relevant results from \cite[Section~2.3]{CI}.

\begin{thm}  \label{tnCIpaper}
Let $B$ be a full matrix algebra over $B_{e,ne}$.
\begin{enumerate}
 \item $B$ is a terminal $k\llbracket u,v \rrbracket$-order and, in particular, regular.
 \item The ramification curves of $B$ are the lines $u=0$ and $v=0$ and the corresponding ramification indices are $e,ne$.
 \item Let $\eps_i \in B_{e,ne}$ be a diagonal primitive idempotent in the matrix form above. Then $B_{e,ne}\eps_i$ is an indecomposable projective $B_{e,ne}$-module and, up to isomorphism, there are no others. 
 \item $B_{e,e}$ is a maximal order and represents an order $e$ element of the Brauer group $H^2_{et}(\spec k\llbracket u,v \rrbracket,k((u,v))^*)$. 
\end{enumerate}
Any terminal order (in the complete local case) is a full matrix algebra over some $B_{e,ne}$.
\end{thm}
Note that when $e=1$ and $n>1$ the orders $B_{e,ne}$ are ramified over a single line.  We also recall that surfaces have terminal singularities if and only if they are smooth \cite[Thm.~4.5(1)]{KM}.  We further note that having terminal ramification data implies that a reflexive order has global dimension two and is a full matrix algebra over $B_{e,ne}$~\cite[Thm.~2.12]{CI}.  For these reasons,
we interpret terminal $k\llbracket u,v \rrbracket$-order as noncommutative analogues of complete local smooth surfaces $\spec k\llbracket u,v \rrbracket.$ This motivates the following definition.
\begin{defn} \label{dqcurve}
A {\em noncommutative plane curve} is an algebra of the form $B/(f)$ where $B$ is a terminal $k\llbracket u,v \rrbracket$-order and $f \in Z(B)$. 
\end{defn}

Our discussion leads to the following

\begin{prop}  \label{pFRTcurves}
 Let $A = B_{e,ne}/(f)$ be a noncommutative plane curve. Suppose that it satisfies Hypothesis~\ref{aram} (for example this occurs when $A$ has finite representation type). Then $\Lam:= B[\sqrt{f}]*G_f$ is a normal order with ramification indices $e$ on $C_u:u=0$, $ne$ on $C_v:v=0$ and 2 on the components $C_1,\ldots, C_s$ of $f=0$. Let 
$$\Delta = (1-\tfrac{1}{e})C_u + (1-\tfrac{1}{ne})C_v + \sum_{i=1}^s\tfrac{1}{2}C_i $$
Then $A$ has finite representation type if and only if $(\spec k\llbracket u,v \rrbracket, \Delta)$ is a log terminal surface as classified in Theorem~\ref{tlogterm}. 
\end{prop}

\section{Noncommutative plane curves} \label{scurves} 

In \cite{GK}, Greuel-Kn\"{o}rrer showed that the plane curves of finite representation type (in the complete local case) are precisely the simple singularities.  All complete local commutative orders of finite representation type are classified in \cite{J,DR}, and these are known to lie between the simple plane curve singularities and their integral closures.
In this section, we wish to extend this result by first classifying noncommutative plane curves of finite representation type. We continue by computing their AR-quivers as subquivers of corresponding log terminal orders.

We note that local orders over commutative complete discrete valuation rings of finite representation type are classified in \cite{HN} building on work of \cite{DK}.  The Auslander-Reiten quivers of Gorenstein orders of finite representation type are classified in \cite{W}, which has been improved in~\cite{Luo}.
There are also interesting classifications of B\"ackstr\"om orders of finite representation type and their quivers in~\cite{RiRo},
and tiled orders due to Zavadskii-Kirchenko~\cite{ZaKi}, and summarized in~\cite[Ch.13.2,Thm.13.9]{Simson}.
We do not address how our classification interacts with the classifications cited above.  It would be interesting to study this.

We continue the notation from the previous section. Let $B$ be a terminal $R$-order where $R = k\llbracket u,v \rrbracket$. Without loss of generality, we may assume that $B = B_{e,ne}$ since we are only interested in facts up to Morita equivalence. We let $0 \neq f \in R, R[\sqrt{f}] = R[z]/(z^2 - f),B[\sqrt{f}]=B[z]/(z^2-f) = R[\sqrt{f}] \otimes_R B$ and $G_f = \langle \s \rangle$, the Galois group of $R[\sqrt{f}]/R$. Recall $\eps = \frac{1}{2}(1-\sigma)$.
\begin{prop}  \label{pnodel} 
Let $\Lam = B[\sqrt{f}] * G_f$. The number of non-isomorphic indecomposable summands of $\Lambda(1-\eps)$ is $n$. 
\end{prop}
\noindent\textbf{Proof.} 
Let $\eps_1,\ldots,\eps_n\in B$ be the diagonal primitive idempotents as in Theorem~\ref{tnCIpaper}. Then the indecomposable summands of $\Lambda(1-\eps)$ are $\Lambda(1-\eps)\eps_i$ where $i=1,\ldots ,n$ since $\eps\Lambda(1-\eps)\eps_i, (1-\eps)\Lambda(1-\eps)\eps_i$ are indecomposable over $B$. $\qed$

Together with Proposition~\ref{pmfact}, the proof of this proposition essentially tells us that the AR-quiver of the noncommutative plane curve $B/(f)$ is the AR-quiver of the log terminal order $B[\sqrt{f}] * G_f$ with the AR-quiver of the terminal order $B$ deleted. Now the AR-quiver of the log terminal order can be computed as a McKay quiver and the key question now is, what AR-subquiver corresponds to $B$?

To answer this question, we need to use the theory of Artin covers as introduced in \cite{A86}. We refer the reader unfamiliar with this concept, to the exposition given in \cite[Section~3]{CHI}. Suffice for now to note that given an order $\Lam$ with $R = Z(\Lam)$ and ramified Galois cover $S/R$ whose ramification indices $r_i$ divide the corresponding ramification indices $e_i$ of $\Lam$, there is a naturally associated $S$-order $\Gamma$, called the {\em Artin cover} such that the following hold:
\begin{enumerate}
 \item There is an extension of the $G$ action from $S$ to $\Gamma$ such that $\Gamma^G = \Lam$.
 \item The ramification indices of $\Gamma$ are $e_i/r_i$. In particular, $\Gamma$ is Azumaya in codimension one if all the $r_i = e_i$. 
\end{enumerate}

The Artin cover of $R[\sqrt{f}]* G_f$ with respect to $R[\sqrt{f}]/R$ is $\Gamma_f = R[\sqrt{f}]\otimes_k k^{2 \times 2}$ and $R[\sqrt{f}]* G_f = \Gamma_f^{G_f}$ where the action of $G_f$ on $\Gamma_f$ can be described as follows. Let $c: G_f \lm \GL_2$ be the representation which sends the non-trivial element $\s \in G_f$ to $c_{\s} = \left(\begin{smallmatrix} 1 & 0 \\ 0 & -1 \end{smallmatrix}  \right)$. Then $\s$ acts by the Galois action on $R[\sqrt{f}]$ and conjugation by $c_{\s}$ on $k^{2 \times 2}$.  

Let $G_B = \Z/e\Z \times \Z/ne\Z$  and $S_B$ be the $G_B$-cover of $R$ that has the same ramification data as $B$, that is, it has ramification index $e$ on one coordinate axis and ramification index $ne$ on the other. The Artin cover of $B$ with respect to $S_B/R$ is $\Gamma_B = S_B\otimes_k k^{ne \times ne}$. As in the previous case, $B = \Gamma_B^{G_B}$ for some action $G_B$ given by a cohomology class $\beta \in H^1(G_B, \PGL_{ne})$ which can be described as follows. Recall that $B$ is determined by a primitive $e$-th root of unity $\zeta$ which corresponds precisely to a generator $d\beta \in H^2(G_B,\mu) \simeq \mu_e$. This lifts to our desired cohomology class $\beta \in H^1(G_B, \PGL_{ne})$. Furthermore, Kummer theory tells us that there exists a central extension $G'_B$ of $G_B$ by $\mu_e$, such that $\beta$ lifts to an element of $b \in H^1(G'_B,\GL_{ne})$, that is, an honest group representation of $G'_B$. 

Now $S_B \otimes_R R[\sqrt{f}]$ is a normal domain since the ramification curves of $S_B$ and $R[\sqrt{f}]$ are distinct. From Corollaries~\ref{cFRT1} and \ref{cFRT2} below, one checks that it is actually a canonical singularity (we will perform this later) so has a Galois cover $S$ which is smooth and the corresponding Galois group $N$ is a subgroup of $\SL_2$. Now $S/R$ is also Galois and the corresponding Galois group $G$ is an extension of $G_B \times G_f$ by $N$. 

\begin{prop}  \label{pArtincover}  
The Artin cover of $\Lambda = B \otimes_R (R[\sqrt{f}]*G_f)$ with respect to $S/R$ is 
$$\Gamma = \Gamma_B \otimes_{S_B} S \otimes_{R[\sqrt{f}]} \Gamma_f = S \otimes_k k^{ne\times ne} \otimes k^{2 \times 2}.$$ 
The action of $G$ on $\Lambda$ is the natural one, that is, $G$ acts via the Galois action on $S$ and via $G_B, G_f$ on $\Gamma_B,\Gamma_f$ respectively. 
\end{prop}

Consider the pullback diagram of group extensions
$$\diagram 
1 \rto & \mu_e \rto \ddouble & G' \rto \dto & G \rto \dto & 1 \\
1 \rto & \mu_e \rto & G_B' \rto & G_B \rto & 1 
\enddiagram$$
where the bottom row is given by $d\beta$. Below, we abuse notation and let $b,c$ denote the representations of $G'$ obtained by composing $b,c$ defined above, with the natural homomorphisms $G' \lm G_B$ and $G' \lm G_f$. 

\begin{thm}  \label{tAR}  
The AR-quiver of $A = B/(f)$ is obtained from $\McK_{d\beta}(G)$ by deleting the vertices where $\ker :=\ker(G' \lm G_B')$ acts trivially. 
\end{thm}
\noindent\textbf{Proof.} In view of Proposition~\ref{pmfact}, we need only use the identification of the AR-quiver of $\Lambda = B[\sqrt{f}] * G_f$ with $\McK_{d\beta}(G)$ as found in \cite[\S4]{CHI}, and track which vertices $\Lambda(1-\eps)$ correspond to. First note that the number of indecomposable projective $B$-modules is $n$ by Theorem~\ref{tnCIpaper}iii) so this is also the number of vertices in $\McK_{d\beta}(G_B)$. This in turn, is also the number of vertices of $\McK_{d\beta}(G)$ where $\ker:=\ker(G' \lm G_B')$ acts trivially and, by Proposition~\ref{pnodel} the number of indecomposable summands of $\Lam (1-\eps)$. It thus suffices to show that $\ker$ acts trivially on the representation of $G'$ corresponding to $\Lam(1-\eps)$. 

This relies on the Artin cover $\Gamma$ of $\Lambda$ with respect to $S/R$ described in Proposition~\ref{pArtincover}. We first use the reflexive Morita equivalence 
$$ \Hom_{\Lambda}(\Gamma,-): \MCM(\Lambda) \xrightarrow{\sim} \MCM(\Gamma*G) .$$
Now one verifies using \cite[Theorem~2.15]{A86} that the trace map induces an isomorphism $\Gamma \simeq \Hom_{\Lambda}(\Gamma,\Lambda)$ of $(\Gamma*G,\Lambda)$-bimodules. We wish thus to see which vertices correspond to the $\Gamma*G$-module $\Gamma(1-\eps)$. 

Let $\delta \in k\mu_e$ be the central idempotent corresponding to the character $\mu \hookrightarrow k^*$. Consider a degenerate action of $G'$ on $\Gamma = S \otimes_k k^{ne \times ne} \otimes_k k^{2 \times 2}$ where $G'$ acts as usual on $S$ but trivially on $k^{ne \times ne} \otimes_k k^{2 \times 2}$. We let $\Gamma \# G'$ be the skew group algebra obtained from this degenerate action. Our next step is to use the algebra isomorphism $\phi:\delta (\Gamma \# G') \xrightarrow{\sim} \Gamma* G: \delta g \mapsto b_g^{-1} g$. We may restrict via $\phi$ to obtain the $\delta (\Gamma \# G')$-module $\Gamma(1-\eps)_{\phi}$. 

The next step is to use the Morita equivalence
$$ \Gamma \# G'-\Mod \simeq S * G': M \mapsto (1 \otimes \eps_B \otimes \eps_f)M$$
where $1 \otimes \eps_B \otimes \eps_f \in S \otimes k^{ne \times ne} \otimes k^{2 \times 2} = \Gamma$ is the idempotent obtained from primitive idempotents $\eps_B \in  k^{ne \times ne}, \ \eps_f \in k^{2 \times 2}$. The $G'$-module corresponding to $\Lam(1-\eps)$ is thus 
$$ V = k \otimes_S (1 \otimes \eps_B \otimes \eps_f)\Gamma(1-\eps)_{\phi} = 
k \otimes \eps_B k^{ne \times ne} \otimes \eps_f k^{2 \times 2} (1-\eps).$$
We work through how $G'$ acts on each of the tensor factors above. Of course $G'$ acts trivially on $k$ so the same is true for $\ker$. We now consider how $g \in G'$ acts on the tensor factor $k^{ne \times ne}$. This is given by the action of $\phi(g) = b_g^{-1} \overline{g} \in \Gamma*G$ where $\overline{g}$ is the image of $g$ in $G$. Now $\overline{g}$ acts by conjugation by $b_g$ so $g$ acts by right multiplication by $b_g^{-1}$. Clearly $\ker$ acts trivially on $k^{ne \times ne}$ and hence also on $\eps_Bk^{ne \times ne}$. Finally, we see similarly that $g$ acts on $k^{2 \times 2} (1-\eps)$ by right multiplication by $c_g^{-1}$. This is trivial since $1-\eps$ is the idempotent corresponding to the trivial character of $G_f$. $\qed$

Theorem~\ref{tlogterm} and Proposition~\ref{pmfact} immediately gives a classification result for noncommutative plane curves of finite representation type. The result is given in the next two corollaries. We will also record the group $G$ and $\ker = \ker(G \lm G_B) \simeq \ker(G' \lm G'_B)$ of Theorem~\ref{tAR} from which one can compute the AR-quivers. The notation for groups is as in Section~\ref{sCoxclass}.


\begin{cor}  \label{cFRT1}  
Let $B$ be a terminal $k\llbracket u,v \rrbracket$-order with ramification indices $e>1,ne$ along the coordinate axes $u=0,v=0$ respectively. Let $0 \neq f \in k\llbracket u,v \rrbracket$ and $C_f$ be the curve $f=0$. Then $B/(f)$ has finite representation type if and only if one of the following occurs.
\begin{enumerate}
\item $e=3,n=1$ and $C_f$ is a smooth curve which is not tangential to either coordinate axis. $G = ST7 = (12,12;E_6,E_6)$ and $\ker = (4,4;D_4,D_4)$. 
\item $e=2$ and either $C_f$ is not tangential to $v=0$ or $n=1$. If $C_f$ intersects $u=0$ with multiplicity $r$ then $G = G(4r,2r,2)$ and $\ker = G(2r,2r,2)$. 
\end{enumerate}
\end{cor}

The case where $B$ is generically split is given in the next corollary.
\begin{cor} \label{cFRT2}  
Let $B$ be a terminal $k\llbracket u,v \rrbracket$-order with ramification index $n>1$ along the coordinate axis $L:v=0$. Let $0 \neq f \in k\llbracket u,v \rrbracket$ and $C_f$ the curve $f=0$. Then $B/(f)$ has finite representation type if and only if one of the following occurs.
\begin{enumerate}
\item $C_f$ is a smooth curve intersecting $L$ with multiplicity $r$ and one of the following occurs: 
 \begin{enumerate}
 \item $n=2$. Then $G = G(2r,2r,2) = (4,2;D_{r+2},A_{2r-1})$ and $\ker = G(r,r,2)$.
 \item $n=3,r=3$. Then $G = ST6 = (12,4;E_6,D_4)$ and $\ker = (4,4;D_4,D_4)$.
 \item $n=3,r=4$. Then $G = ST14 = (12,6;E_7,E_6)$ and $\ker = (4,2;E_7,E_6)$.
 \item $n=3,r=5$. Then $G = ST21 = (12,12;E_8,E_8)$ and $\ker = (4,4;E_8,E_8)$.
 \item $n=4, r = 3$.  Then $G = ST9 = (8,8;E_7,E_7)$ and  $\ker = (8,4;E_7,E_6)$. 
 \item $n=5, r=3$. Then $G = ST17 = (20,20;E_8,E_8)$ and $\ker = (4,4;E_8,E_8)$.
 \item $r = 2$. Then $G = G(n,1,2)$ and $\ker = G(n,n,2)$. 
 \item $r=1$. Then $G = \mu_n \times \mu_2$ and $\ker = \mu_2$.  
 \end{enumerate}
\item $C_f$ is a cusp of type $A_{2r}$ which is not tangential to $L$. Then $G = G((2r+1)n,2r+1,2)$ and $\ker = G((2r+1)n,(2r+1)n,2)$
\item $C_f$ is the union of two non-tangential smooth curves and if $r= C_f.L - 1$ then one of the following occurs: 
 \begin{enumerate}
 \item  $n=2$. Then $G = G(4r,2r,2)$ and $\ker = G(2r,r,2)$. 
 \item $n=3, r= 2$.  Then $G = ST15 = (12,12;E_7,E_7)$ and $\ker = (4,4;E_7,E_7)$.
 \item $r=1$. Then $G = G(2n,2,2)$ and $\ker = G(2n,2n,2)$. 
 \end{enumerate}
\item $C_f$ is a type $A_{2r-1}$ node which is not tangential to $L$. Then $G = G(2rn,2r,2)$ and $\ker = G(2rn,2rn,2)$.
\item  $C_f$ is an ordinary cusp which is tangential to $L$ and $n=2$. Then $G = ST13 = (4,4;E_7,E_7)$ and $\ker = ST12 = (4,2;E_7,E_6)$.
\end{enumerate}
\end{cor}

The case where $B$ is ramified on $L:v=0$ only and $C_f$ is a smooth curve intersecting $L$ with multiplicity two was shown to be finite representation type in \cite{CC}. Indeed, in that case, the indecomposable Cohen-Macaulay modules were explicitly computed without use of the McKay quiver.

\section{Appendix}  \label{sappendix}
We present some of the McKay graphs of finite subgroups of $\GL_2$.

\begin{tabular}{c|c|c|c}
{\includegraphics[scale=0.2]{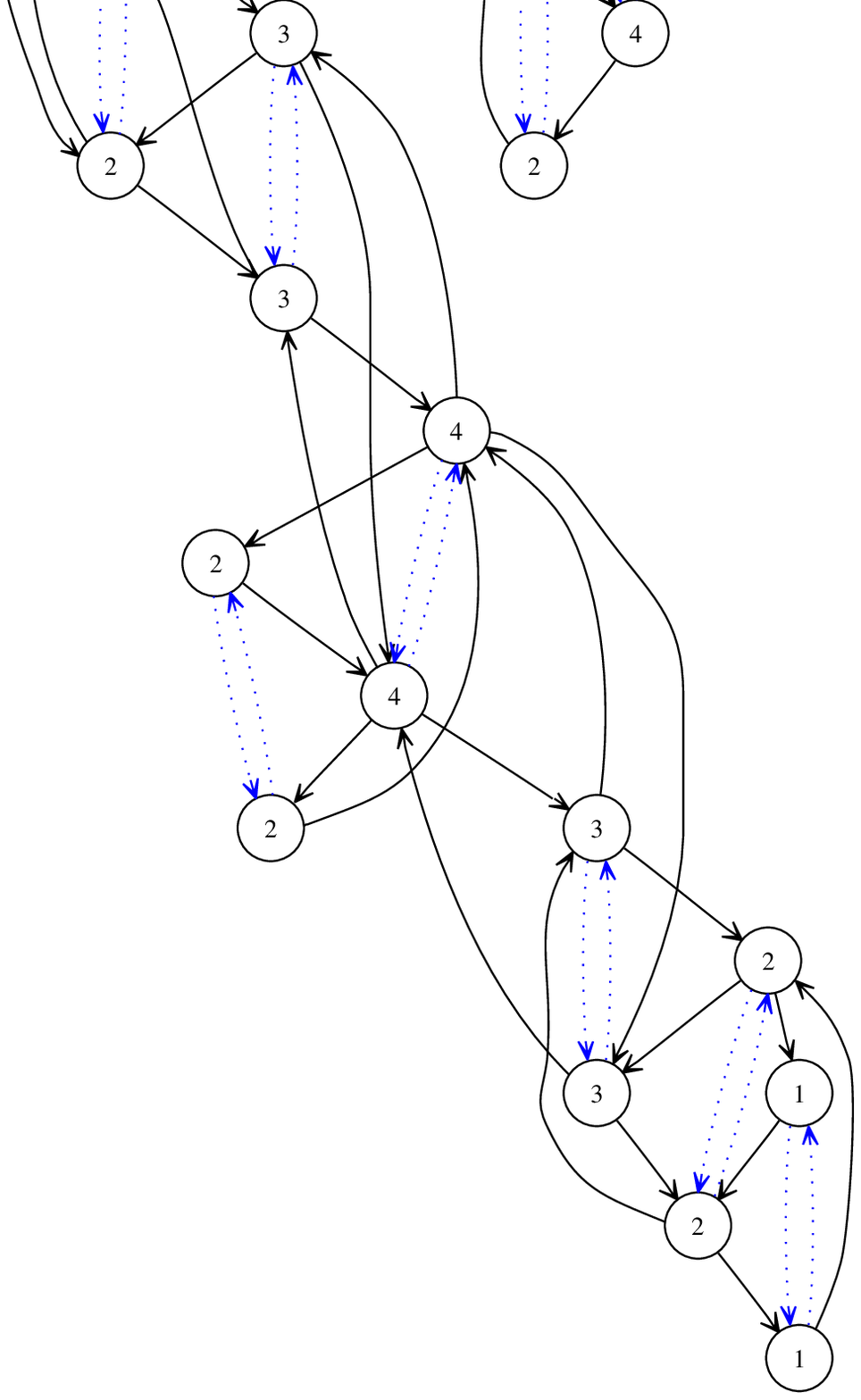}} & {\includegraphics[scale=0.3]{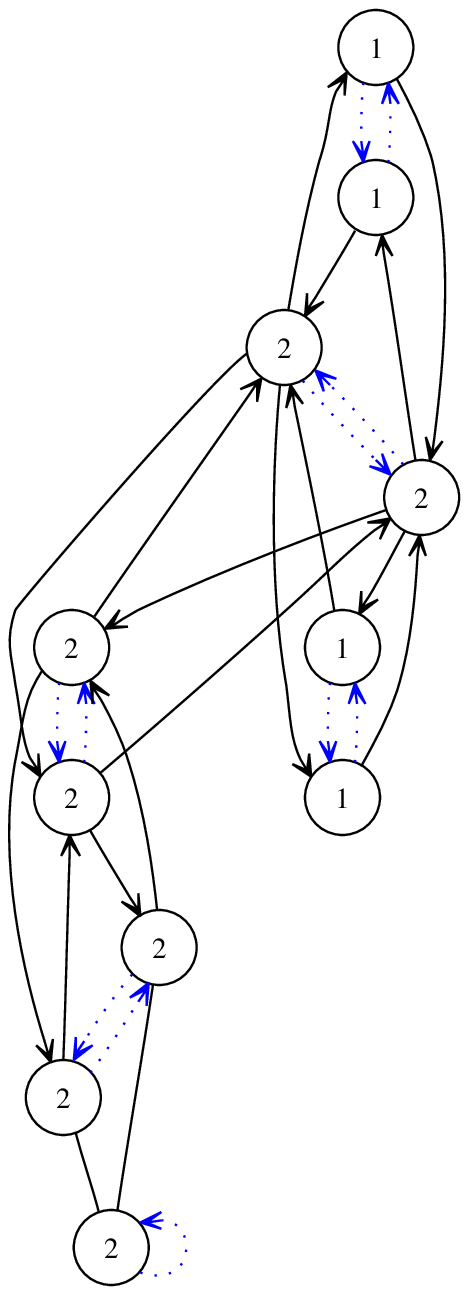} } & {\includegraphics[scale=0.3]{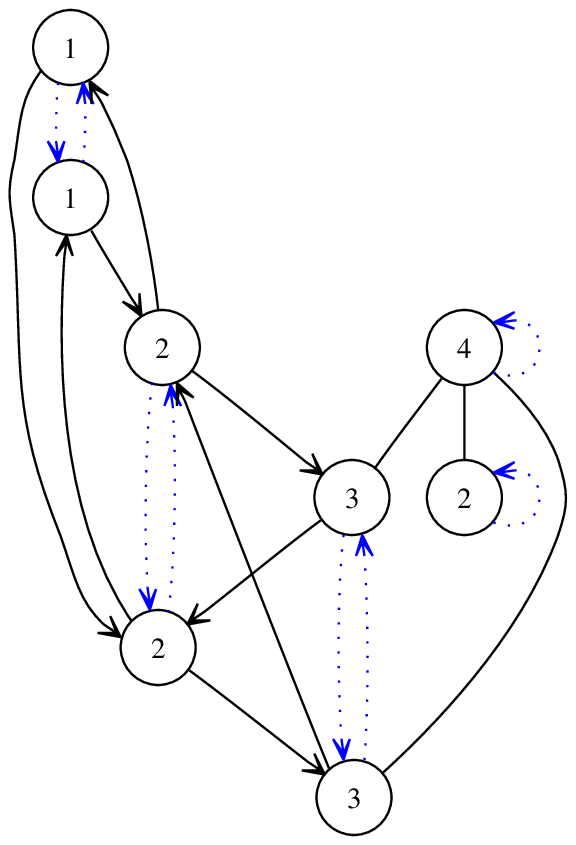}
} & \includegraphics[scale=0.3]{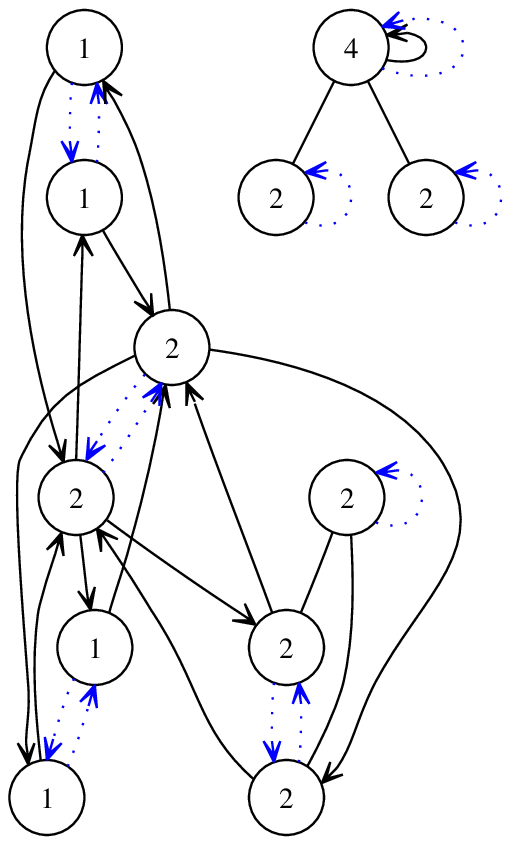}\\
$E_7^2$ and $F_{42}$ & $CD_5$ & $F_{41}$ & $CD_4$ and $DT_2$
\end{tabular}

\begin{tabular}{c|c}
 {\includegraphics[scale=0.2]{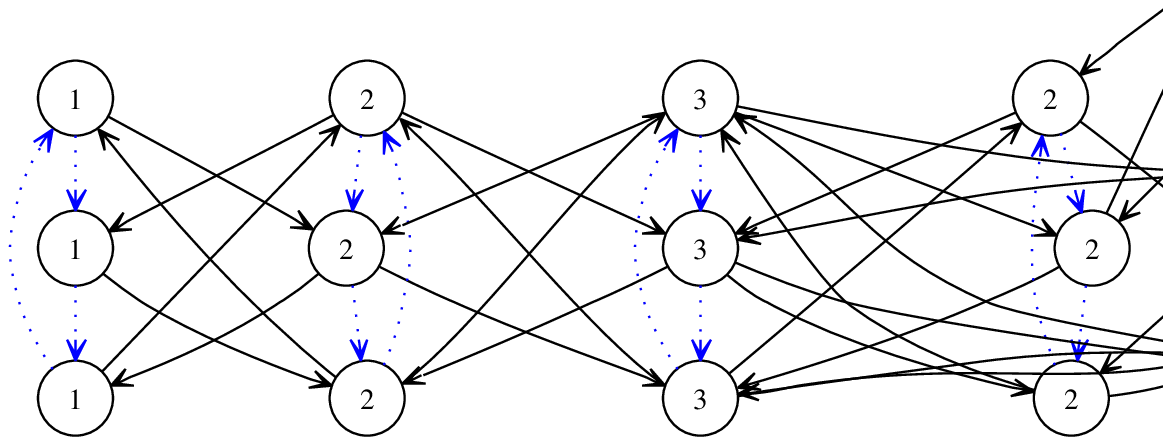}} & {{\includegraphics[scale=0.2]{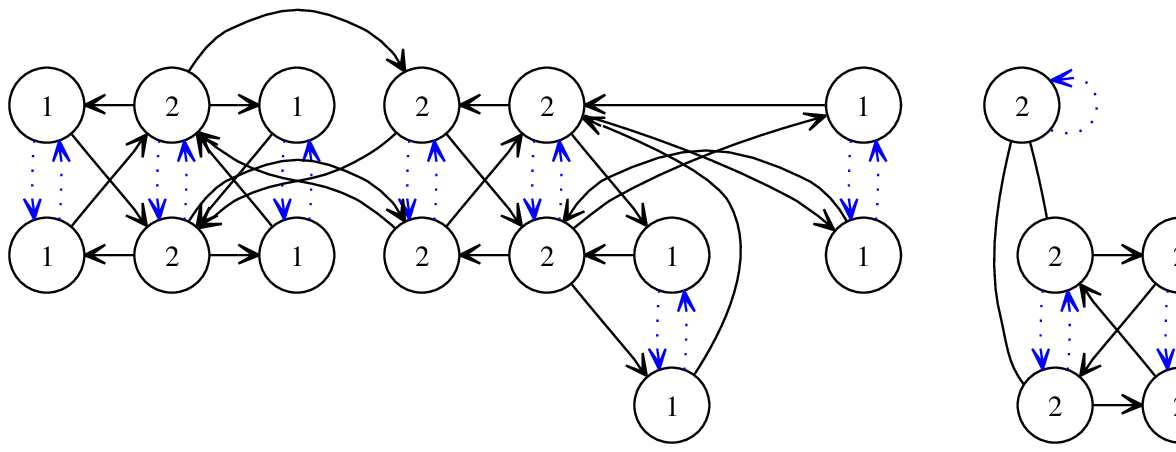}}}\\
 {$E_6^3$ and $G_{22}$}  & {$D_6^2$ and $C_4$ and $BD_3$ and $BD_3$}
\end{tabular}

\begin{tabular}{c|c}
 {{\includegraphics[scale=0.2]{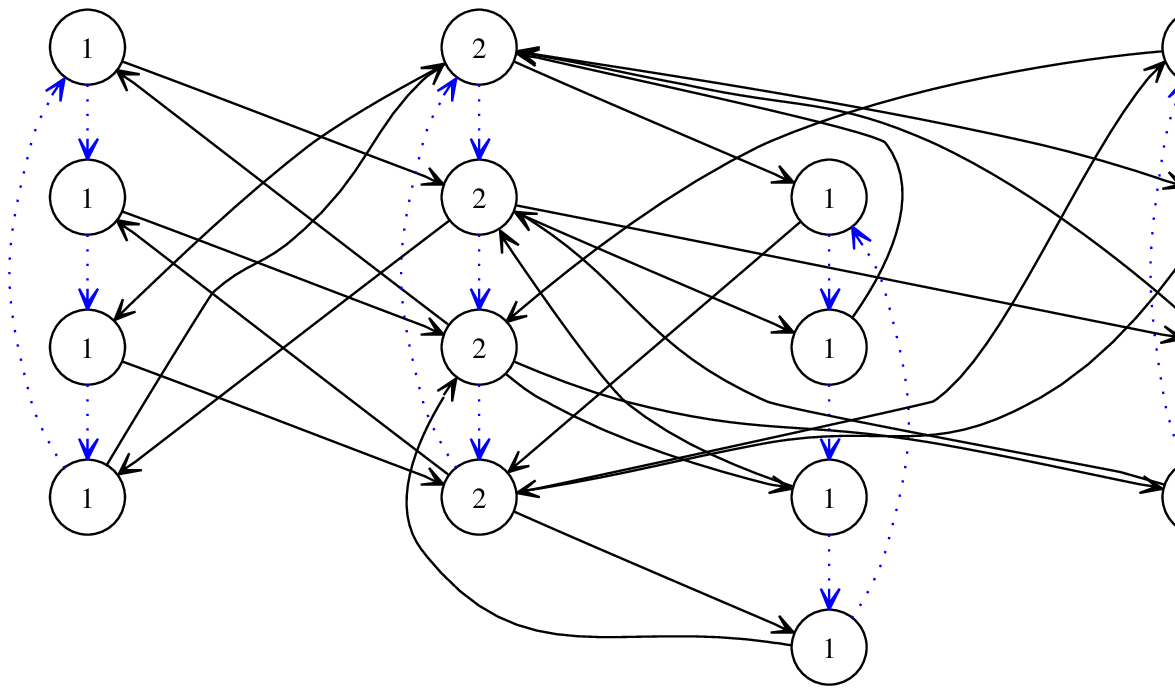}}}  & {{\includegraphics[scale=0.3]{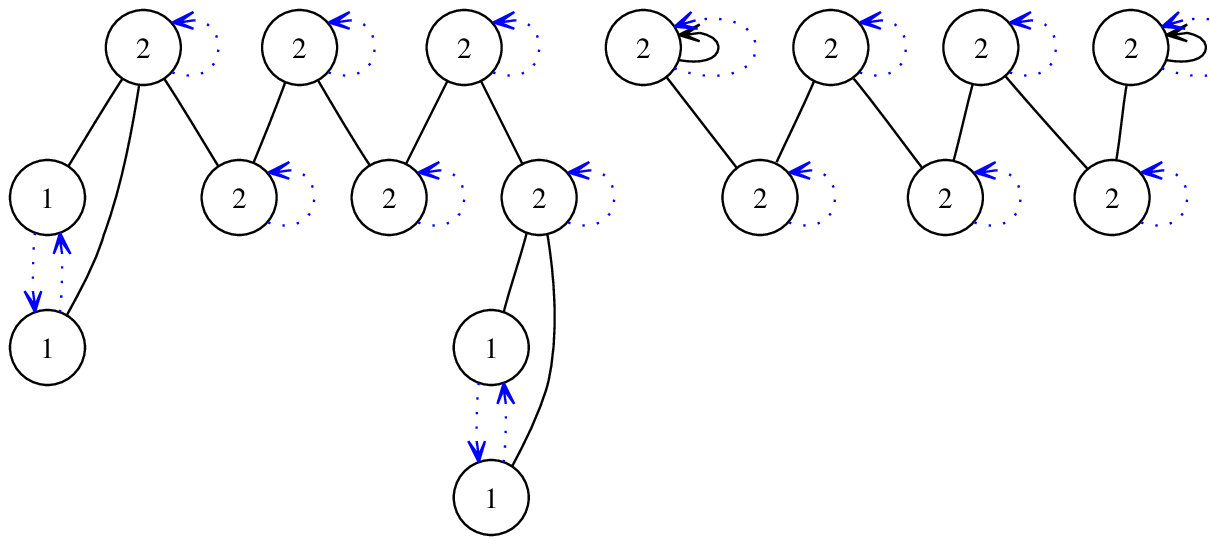}}} \\
{$CD^5_2$ and $BC_2$} & {$B_7$ and $T_6$}
\end{tabular}

\begin{tabular}{c|c|c}
  {{\includegraphics[scale=0.3]{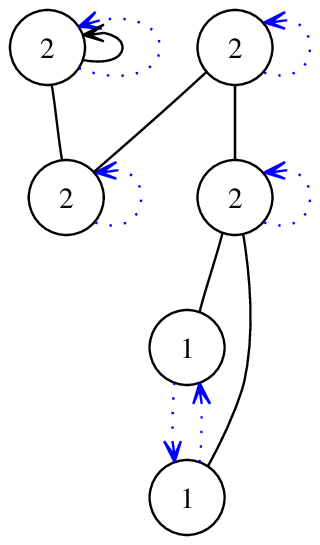}}} & {{\includegraphics[scale=0.3]{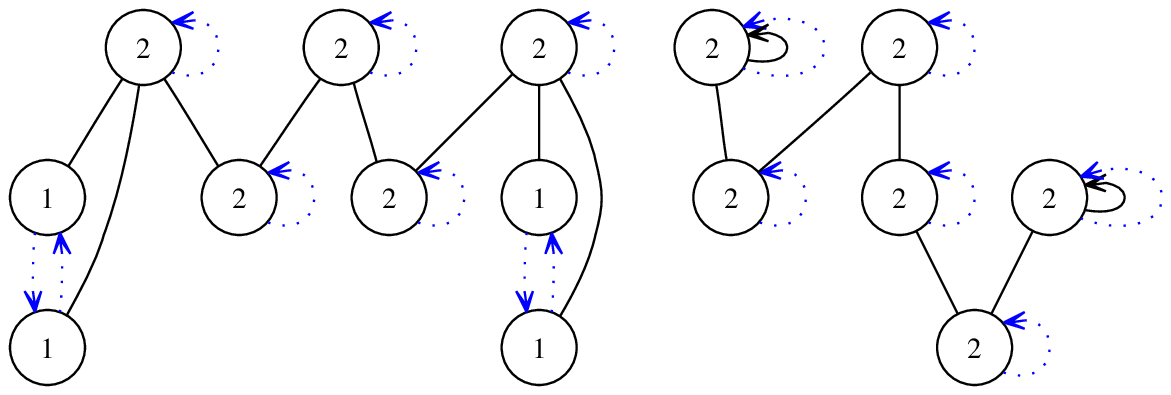}}}& {\includegraphics[scale=0.3]{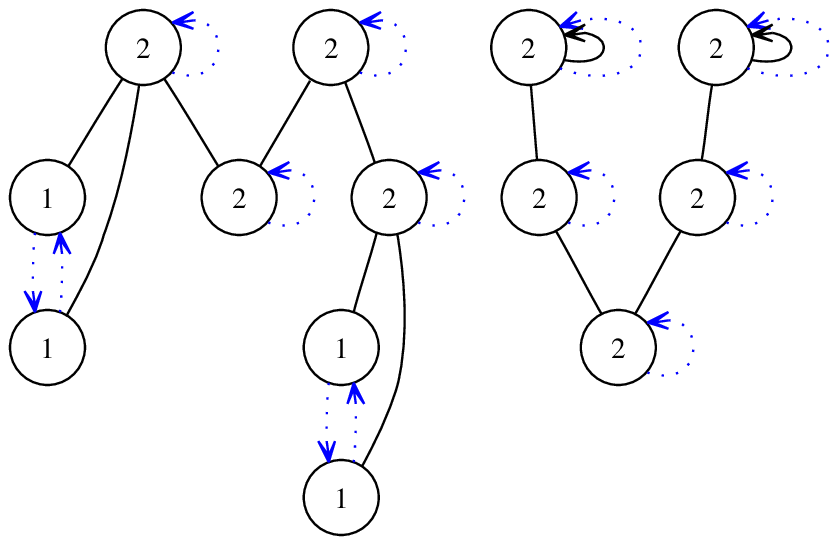}}\\
$BT_4$ &   $B_6$ and $T_5$ & {$B_5$ and $T_4$}
\end{tabular}

\begin{tabular}{c|c}
{{\includegraphics[scale=0.3]{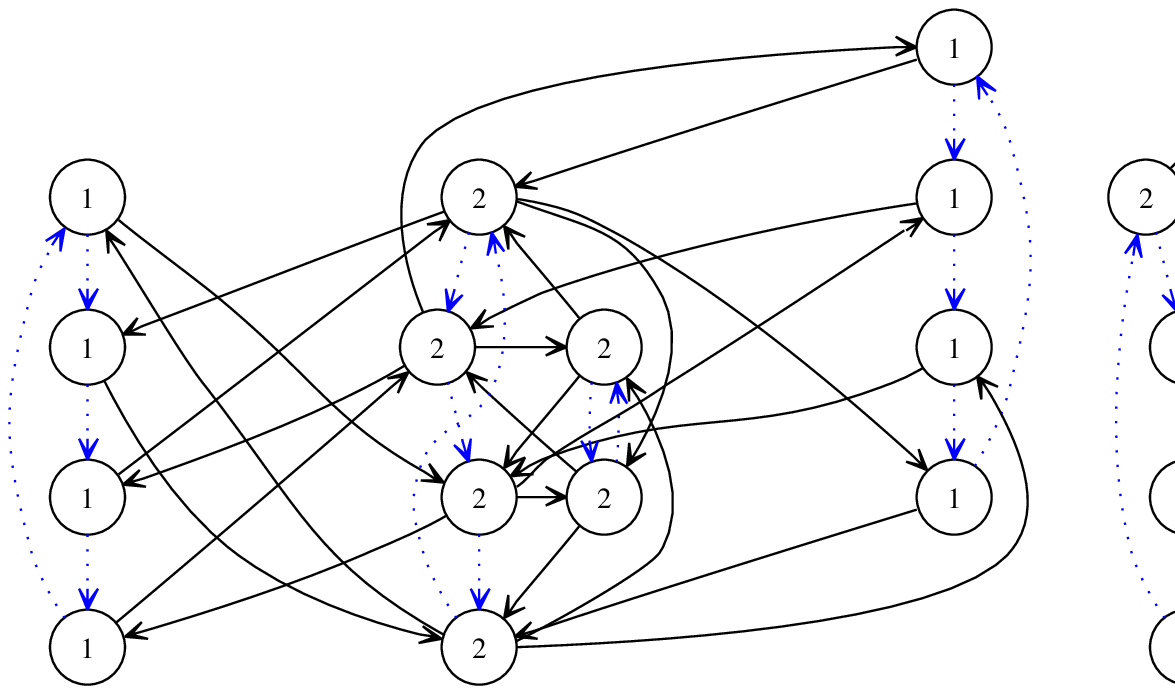}}} & {{\includegraphics[scale=0.3]{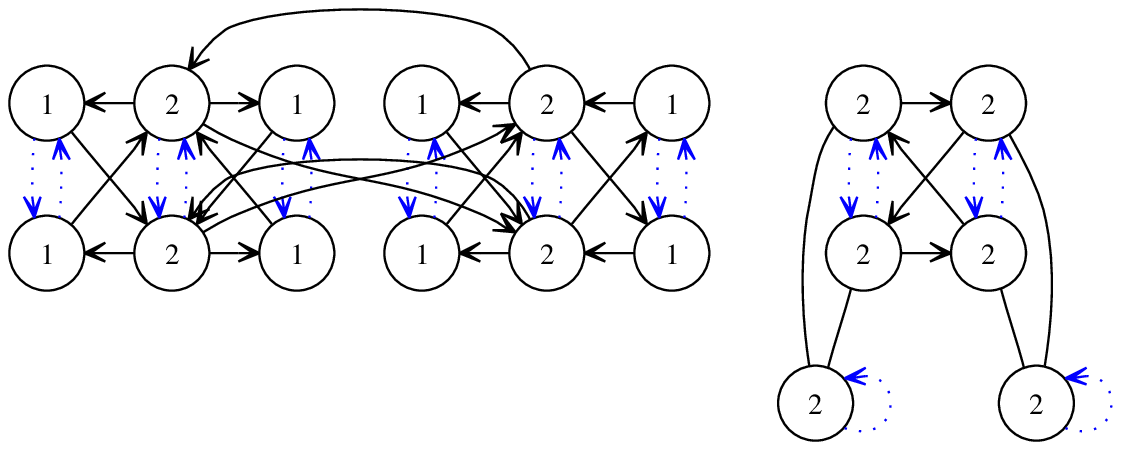}}}\\
{$CD^2_3$ and $A_{11}$} & {$C_3$}
\end{tabular}

In these tables we use the following abbreviations. A positive integer we abbreviate $\Z/a$ as simply $a$ and similarly we write $a \oplus b$ for the group $\mathbb{Z}/a \oplus \Z/b$.  As before $a \vee b$ for the greatest common divisor of $a$ and $b$.
Below is a list of the cohomology groups $H^1(G,k^*), H^2(G,k^*)$ for all subgroups $G$ of $\GL_2$ as computed in Sections~\ref{sARsplit} and \ref{sARgeneral}. 


$$ \begin{array}{l|l|l}
\text{type}     & H^1 & H^2 \\
\hline
B_n^m    &  \begin{array}{ll}
	   4m & n \mbox{ odd }, m \mbox{ even } \\ 
	   2 \oplus 2m & \mbox{ else }  \\ 
	   \end{array}     & (m+1) \vee 2 \\
\hline
D_n^m    & \begin{array}{ll}
           2 \oplus 2 \oplus m & n \mbox{ even } \\ 
	   4 \oplus m & n \mbox{ odd }, m \mbox{ odd } \\ 
	   2 \oplus 2m & n \mbox{ odd }, m \mbox{ even } 
	   \end{array}  
          & (m \vee 2) \oplus (m \vee n \vee 2) \\
\hline E_6^m    & 3 \oplus m   & (m \vee 3) \\
\hline E_7^m    & 2 \oplus m   & (m \vee 2)\\
\hline E_8^m    & m  & 1   \\
\hline CD_n^m   & 2 \oplus 2m & (m+n-1) \vee 2 \\
\hline F_{41}^m & 2m &  1  \\ 
\hline G_{21}^m & 3m &  1  \\ 
\hline BT_n^m   & 2(2m-1) &  1 
\end{array}$$  


The next table lists the information about the pseudo-reflection subgroup of $G = G(Z_1,Z_2;G_1,G_2)$ as mentioned in Theorem~\ref{tramofklt}. In the cases where there are two parameters $m,n$ we let $\overline{m} = m/(m\vee n), \overline{n} = n/(m\vee n)$. The notation $\frac{1}{a}(b,c)$ denotes the cyclic group 
$ \left\langle 
\left(\begin{smallmatrix}
 \zeta_a^b & 0 \\ 0 & \zeta_a^c
\end{smallmatrix}\right)
\right\rangle$.
 $$ \begin{array}{l|l|l}
 \mbox{ type } & RG & G/RG \\
 \hline
 A^c_{m,n} & m \oplus n\vee c & n/n\vee c \\
B^m_n &    \begin{array}{ll}
         G(2n,2\overline{n},2) & m \vee 2 = 1 \\
         \mu^2_{m \vee n} & m \vee 2 = 2 \\
                          & 
         \end{array}  & 
         \begin{array}{ll}
         \frac{1}{\overline{m}}(\overline{n},1) & m \vee 2 = 1 \\
         (2\overline{m},2\overline{m};D_{\overline{n}+2},D_{\overline{n}+2}) & m \vee 2 = 2, \overline{m} \vee 2 = 1 \\
         (4\overline{m},2\overline{m};D_{\overline{n}+2},A_{2\overline{n}}) & m \vee 2 = 2, \overline{m} \vee 2 = 2\\
         \end{array}      \\ \hline
D_{n+2}^m      &   \begin{array}{ll}
         G(2n,2\overline{n},2)  & m \vee 2 = 2 \\
         \mu^2_{m \vee n} & m \vee 2 = 1
         \end{array}  & 
         \begin{array}{ll}
         \frac{1}{\overline{m}}(\overline{n},1)  & m \vee 2 = 2 \\
          (2\overline{m},2\overline{m};D_{\overline{n}+2},D_{\overline{n}+2}) & m \vee 2 = 1
         \end{array}    \\ \hline
E_6^m  &   \begin{array}{ll}
        1 & m \vee 6=1 \\
       D_4^2 & m \vee 6=2 \\
       E_6^{m \vee 6} & \mbox{else} \\
         \end{array}   &  
      \begin{array}{ll}
             E_6^m & m \vee 6=1 \\
           3m/2  & m \vee 6=2 \\
       m/m \vee 6 & \mbox{else} \\
    \end{array}        \\ \hline
E_7^m  &  \begin{array}{ll}
        1 & m \vee 12=1 \\
       E_6^3 & m \vee 12=3 \\
      E_7^{m \vee 12} &\mbox{else} \\
         \end{array}   &  
      \begin{array}{ll}
             E_7^m & m \vee 12=1 \\
            2m/3 & m \vee 12=3 \\
       m/m \vee 12 & \mbox{else} \\
    \end{array} 
    \\ \hline
E_8^m  &  \begin{array}{ll}
        1 & m \vee 30=1 \\
       E_8^{m\vee 30} & m \vee 30 
      \neq 1 \end{array}  &  
      \begin{array}{ll}
            E_8^m & m \vee 30=1 \\
           m/m \vee 30 & m \vee 30
           \neq 1 \end{array} 
           \\ \hline
CD^m_{n+1}  &  \begin{array}{ll}
              G(2n,2\overline{n},2) & m \vee 2 = 2, (\overline{m}\overline{n}) \vee 2 = 2\\
              G(2n,\overline{n},2) & m \vee 2 = 2, (\overline{m}\overline{n}) \vee 2 = 1 \\
              G(2n,2\overline{n},2) & m \vee 2 = 1, (\overline{m}\overline{n}) \vee 2 = 2 \\
              G(2n,\overline{n},2) & m \vee 2 = 1, (\overline{m}\overline{n}) \vee 2 = 1 
             \end{array} & 
             \begin{array}{ll}
              \frac{1}{2\overline{m}}(\overline{m} + \overline{n},1)& m \vee 2 = 2, (\overline{m}\overline{n}) \vee 2 = 2 \\
              \frac{1}{\overline{m}}(\frac{\overline{m} + \overline{n}}{2},1)& m \vee 2 = 2, (\overline{m}\overline{n}) \vee 2 = 1 \\
              \frac{1}{2\overline{m}}(\overline{m} + \overline{n},1) & m \vee 2 = 1, (\overline{m}\overline{n}) \vee 2 = 2 \\
              \frac{1}{\overline{m}}(\frac{\overline{m} + \overline{n}}{2},1)& m \vee 2 = 1, (\overline{m}\overline{n}) \vee 2 = 1   
             \end{array}   \\ \hline
F^m_{41}  &   \begin{array}{ll}
        D_4^2 & m \vee 12=4 \\
       E_6^6 & m \vee 12=12 \\
F_{41}^{m\vee 12} &\mbox{else} \\
         \end{array}   &  
      \begin{array}{ll}
             3m & m \vee 12=4 \\
           m/3  & m \vee 12=12 \\
       m/m \vee 12 & \mbox{else} \\
    \end{array} 
    \\\hline
G^m_{21} &   \begin{array}{ll}
        1 & m \vee 6=3 \\
       D_4^2 & m \vee 6=6 \\
       G_{21}^{m \vee 6} & \mbox{else} \\
         \end{array}   &  
      \begin{array}{ll}
             G_{21}^m & m \vee 6=3 \\
           3m/2  & m \vee 6=6 \\
       m/m \vee 6 & \mbox{else} \\
    \end{array}      \\ \hline
BT^{m}_{\frac{n-1}{2}}  &
                \begin{array}{ll}
                 G(n,\overline{n},2) \\ 
                \end{array}  &
                \begin{array}{ll}
                 \frac{1}{\overline{m}}(\frac{3m+1}{4}\overline{n}, \frac{1+\overline{m}}{2}) & m \equiv 1 \mod 4 \\
                  \frac{1}{\overline{m}}(\frac{m+1}{4}\overline{n}, \frac{1+\overline{m}}{2}) & m \equiv 3 \mod 4 
                \end{array}     \\

\end{array}$$  

\bibliographystyle{alpha}
\bibliography{references}

\end{document}